\documentclass[12pt]{elsarticle}

\usepackage{amssymb}
\usepackage{amsmath}
\journal{Physica D}

\usepackage{xcolor}
\usepackage[colorlinks=true]{hyperref}
\usepackage{multirow}   % for tables
\usepackage{siunitx} %

\usepackage{algorithm2e}
\usepackage{algpseudocode}
\usepackage{graphicx}

\usepackage{amsthm}

\usepackage{graphics}
\usepackage{subcaption}

\usepackage{mathtools}
\mathtoolsset{showonlyrefs}% only print equation numbers if they are referenced (needs mathtools)

\makeatletter
\newcommand*{\defeq}{\mathrel{\rlap{%
                     \raisebox{0.3ex}{$\m@th\cdot$}}%
                     \raisebox{-0.3ex}{$\m@th\cdot$}}%
                     =}
\newcommand*{\eqdef}{=\mathrel{\rlap{%
                     \raisebox{0.3ex}{$\m@th\cdot$}}%
                     \raisebox{-0.3ex}{$\m@th\cdot$}}%
                     }
\makeatother

\newtheorem{lemma}{Lemma}
\newtheorem{problem}{Problem}
\newtheorem{remark}{Remark}

\DeclareMathOperator*{\esssup}{ess\,sup}
\DeclareMathOperator*{\argmax}{arg\,max}

\newcommand{\mathcolon}{\mathbin{:}}

\renewcommand{\H}{%
  \ifmmode%
    H^{\frac{3}{2}-\frac{1}{q}}%
  \else%
    $H^{3/2-1/q}$%
  \fi%
}%
\newcommand{\x}{\ensuremath{\boldsymbol{x}}}
\renewcommand{\d}{\ensuremath{\boldsymbol{d}}}
\renewcommand{\u}{\ensuremath{\boldsymbol{u}}}
\newcommand{\uInit}{\u^{\text{init}}}
\renewcommand{\v}{\ensuremath{\boldsymbol{v}}}
\newcommand{\w}{\ensuremath{\boldsymbol{w}}}
\newcommand{\n}{\ensuremath{\boldsymbol{n}}}

\newcommand\bnabla{\ensuremath{\boldsymbol{\nabla}}}
\newcommand{\0}{\boldsymbol{0}}
\newcommand{\bxi}{\ensuremath{\boldsymbol{\xi}}}
\newcommand{\uInitPDE}{\u^0}
\newcommand{\Bchange}{\bar B}

\usepackage{mathrsfs}
\newcommand{\Retraction}{\mathsf{R}}

\newcommand{\nablaH}{\ensuremath{\nabla^H}}
\newcommand{\nablaL}{\ensuremath{\nabla^L}}
\newcommand{\nablaB}{\ensuremath{\nabla^B}}
\newcommand{\R}{\ensuremath{\mathcal{R}}}
\newcommand{\G}{\ensuremath{\mathcal{G}}}

\newcommand{\largeB}{"large" $B$}
\newcommand{\smallB}{"small" $B$}

\newcommand{\integrandFull}{r}

\newcommand{\TXB}[1][]{\ensuremath{\mathcal{T}_{#1}X_B}}

\def\TukmoXB{\TXB[\u_{k-1}]}
\def\TukXB{\TXB[\u_{k}]}
\def\TuXB{\TXB[\u]}

\newcommand{\hatvp}{\widehat{\v'}}

\newcommand{\hide}[1]{}

\usepackage[color=orange]{todonotes}%leave used. otherwise tikz breaks in rieConjGrad.tex

\colorlet{revOneColor}{violet}

\definecolor{paraviewRed}{rgb}{0.85,0,0}
\definecolor{paraviewGreen}{rgb}{0,0.85,0}
\definecolor{paraviewBlue}{rgb}{0,0,0.85}
\definecolor{paraviewTurquoise}{rgb}{0,0.85,0.85}
\definecolor{paraviewPurple}{rgb}{0.85,0,0.85}
\definecolor{paraviewYellow}{rgb}{0.85,0.85,0}
\definecolor{paraviewOrange}{rgb}{0.85,0.425,0}

\definecolor{matplotlibBlue}{HTML}{1f77b4}
\definecolor{matplotlibOrange}{HTML}{ff7f0e}
\definecolor{matplotlibGreen}{HTML}{2ca02c}
\definecolor{matplotlibRed}{HTML}{d62728}
\definecolor{matplotlibPurple}{HTML}{9467bd}
\def\RR{\mathbb{R}}
\def\TT{\mathbb{T}}
\def\CC{\mathbb{C}}
\def\ZZ{\mathbb{Z}}

\def\E{\mathcal{E}}

\def\S{\mathcal{S}}

\def\tuB{\widetilde\u^B}
\def\tuE{\widetilde{\u}^{\E_0}}
\def\tC{\widetilde{C}}
\def\talpha{\widetilde{\alpha}}

\def\TukXB{\mathcal{T}_{\u_k} X_B}
\newcommand{\hatv}{\ensuremath{\widehat{\v}}}

\begin{document}

\begin{frontmatter}



\title{On the sharpness of bounds on the rate of growth of Lebesgue norms of the velocity in Navier-Stokes flows}


\author[mcmasterAffil]{Fabian Bleitner} %% Author name
\author[mcmasterAffil]{Bartosz Protas} %% Author name

\affiliation[mcmasterAffil]{organization={Department of Mathematics \& Statistics, McMaster University},%Department and Organization
            addressline={1280 Main Street West}, 
            city={Hamilton},
            postcode={L8S 4L8}, 
            state={Ontario},
            country={Canada}}

%% Abstract
\begin{abstract}
    %\tolerance=1000%allows better line breaking
    In this paper we consider solutions $\u$ of the three-dimensional Navier-Stokes system and investigate sharpness of the a priori bound
    \begin{align*}
        \frac{d}{dt}\|\u\|_q^q \leq C\|\u\|_q^{q\frac{q-1}{q-3}}, \qquad q > 3.
    \end{align*}
     This bound is closely related to the Ladyzhenskaya-Prodi-Serrin conditions characterizing classical solutions of the Navier-Stokes system. Velocity fields maximizing the rate of growth $(d/dt)\|\u\|_q^q$ under certain constraints are found as solutions of a suitable optimization problem which is solved numerically using a Riemannian conjugate gradient approach. The results obtained for different $q$ and increasing values of $\|\u\|_q$ indicate that the bound is indeed sharp, up to a numerical prefactor, and therefore cannot be fundamentally improved. 
     Additionally, the results also suggest that the rate of growth $(d/dt)\|\u\|_q^q$ diverges as $q\to 3$.
\end{abstract}



%% Keywords

\end{frontmatter}

%\tableofcontents
%\listoftodos

\section{Introduction}
\label{sectionIntro}
We consider flows of viscous incompressible fluids on a domain without solid boundaries which we assume to be a three-dimensional (3D) torus with unit side length $\Omega = \TT^3 = \RR^3 / \ZZ^3 = [0,1]^3$ equipped with periodic boundary conditions. They are governed by the Navier-Stokes equation given by 
\begin{alignat}{1}
    \partial_t \u + \u\cdot \bnabla \u + \bnabla p - \nu\Delta \u &= \0, \label{nse}
    \\
    \bnabla \cdot \u &= 0, \label{incompressible}
\end{alignat}
where $\u = \u(t,\x)$ is the velocity field, $p = p(t,\x)$ the corresponding pressure and $\nu>0$ is the coefficient of the kinematic viscosity, with $\x \in \Omega$ and $t \ge 0$ denoting, respectively, the space coordinate and the time. Without loss of generality, the constant fluid density is assumed equal to unity. System \eqref{nse}--\eqref{incompressible} is complemented with a suitable initial condition $\uInitPDE$. When there is no risk of confusion, we will use the simplified notation $\u(t) = \u(t,\cdot)$.

System \eqref{nse}--\eqref{incompressible} depends on three physical parameters: the  characteristic length scale and velocity of the flow, and the viscosity coefficient $\nu$. Given the scaling symmetry of system \eqref{nse}--\eqref{incompressible}, they can be combined into a single dimensionless parameter, the Reynolds number. Due to the nature of the problem considered in the present study, we fix $\nu=1$ and use $\|\u\|_q$ as the adjustable parameter.

One of the most fundamental open questions in fluid dynamics is whether the Navier-Stokes equations \eqref{nse}--\eqref{incompressible} always admit globally-defined unique smooth solutions $\u(t)$, $t > 0$, provided the initial data $\uInitPDE$ is sufficiently smooth, or, conversely, if smooth initial data can produce singularities in finite time where system \eqref{nse}--\eqref{incompressible} is no longer satisfied in the classical (pointwise) sense. The Clay Mathematics Institute recognizes this fundamental question as one of its Millennium Problems posed as challenges to the mathematics community \cite{clayMillenniumProblem}.

The study of the problem stated above revolves around the so-called conditional regularity results which are conditions that must be satisfied by the solution $\u(t)$ of \eqref{nse}--\eqref{incompressible} if a singularity does not occur and must be violated if it does. One of the best-known results of this type is the family of the Ladyzhenskaya-Prodi-Serrin conditions \cite{ladyzhenskaya57,prodi59,serrin62} asserting that the solution $\u(t)$ remains smooth on $[0,T]$ as long as
\begin{align}
    \int_0^T \|\u(t)\|_q^p \, dt <\infty, \quad \frac{2}{p}+\frac{3}{q}=1, \qquad q > 3,
    \label{lps}
\end{align}
with the borderline (critical) case $q = 3$ obtained in \cite{Escauriaza2003}
\begin{align}
    \esssup_{0\leq t\leq T}\|\u(t)\|_3 <\infty.
\end{align}
Relation \eqref{lps} implies that if a singularity is to form in a Navier-Stokes flow $\u(t)$ at time $t^\star$, then necessarily 
\begin{align}
    \lim_{T \rightarrow t^\star} \int_0^T \|\u(t)\|_q^p \, dt = \infty
    \label{lpst*}
\end{align}
and the question is then about finding an initial condition $\uInitPDE$ subject to certain constraints such that this can happen. Since whether or not this scenario can occur is a matter of how large the norms $\| \u(t) \|_q$ can grow, it is instructive to consider upper bounds on their rate of growth that were obtained in \cite{robinsonSadowski2014}
\begin{align}
    \frac{1}{q}\frac{d}{dt}\|\u(t)\|_q^q &\leq C\nu^{-\frac{q+3}{q-3}} \|\u(t)\|_q^{q\frac{q-1}{q-3}}, \qquad q > 3,
    \label{estimate}
\end{align}
where the constant $C>0$ depends only on $q$ (hereafter, $C$ will refer to a generic positive constant whose actual value may vary from instance to instance); for completeness, this bound is rederived in \ref{appendixProofOfAprioriBound}. 

As can be verified, cf.~Lemma \ref{lemmaSaturationImpliesBlowUp} in \ref{appendixSaturationImpliesBlowup}, a solution $\u(t)$ of the Navier-Stokes system must saturate bound \eqref{estimate} sufficiently long if a singularity as characterized by \eqref{lpst*} is to occur. Our goal in the present study is to investigate the sharpness of bound \eqref{estimate}. We clarify that this estimate can be declared "sharp" if there exists a constant $C>0$ and a family of fields $\tuB$ parameterized by their $L^q$ norm $B = \| \tuB \|_q$ such that $(d/dt) \| \tuB \|_q^q = C \| \tuB \|_q^{q\frac{q-1}{q-3}}$ as $B \rightarrow \infty$. The exponent on the right-hand side in \eqref{estimate} becomes unbounded as $q \rightarrow 3^+$ and the bound is not applicable when $q = 3$. We are not aware of any bound on $(d/dt) \| \u(t) \|_3$ that can be expressed as a power of  $\| \u(t) \|_3$. While \cite{Tao2020} states a lower bound on how rapidly $\| \u(t) \|_3$ must blow up as a hypothetical singularity is approached, this does not translate to a polynomial estimate on $(d/dt) \| \u(t) \|_3$ analogous to \eqref{estimate}.

Another important conditional regularity result is the enstrophy condition \cite{foiasTemam1989}, which states that the solution $\u(t)$ remains smooth on a time interval $[0,T]$ as long as
\begin{align}
    \sup_{0\leq t\leq T}\|\bnabla \times \u(t)\|_2 < \infty,
    \label{maxEinf}
\end{align}
where the quantity $\|\bnabla \times \u\|_2^2$ is referred to as the enstrophy. In the absence of boundary effects, such as on the three-dimensional torus $\TT^3$ consider here, the equality $\|\bnabla \times \u\|_2^2 = \|\bnabla \u\|_2^2$ holds. Thus, the question whether or not a singularity may form in finite time can again be framed as whether there exists an initial condition $\uInitPDE$ with a finite enstrophy $\|\bnabla \uInitPDE \|_2^2 < \infty$ such that the corresponding Navier-Stokes flow $\u(t)$ violates \eqref{maxEinf} for some $T = t^\star$. As regards how large enstrophy can grow, its rate of change is subject to the bound \cite{luDoering08}
\begin{align}
    \frac{d}{dt}\|\bnabla \u (t)\|_2^2 \leq \frac{1}{2}\left(\frac{3}{\sqrt{2\pi}\nu}\right)^3 \|\bnabla \u(t)\|_2^6.
    \label{jwekrnwer}
\end{align}
Integrating this inequality with respect to time leads to the finite-time estimate
\begin{align}
    %\|\bnabla \u\|_2^2 \leq \frac{\|\bnabla \uInitPDE\|_2^2}{\sqrt{1-C\|\bnabla \uInitPDE\|_2^4t}}.
    \|\bnabla \u(t)\|_2^2 \leq \frac{\|\bnabla \uInitPDE\|_2^2}{\sqrt{1-C\|\bnabla \uInitPDE\|_2^4t}},
    \label{maxE}
\end{align}
where we note that the upper bound on the right-hand side becomes unbounded as $t \rightarrow C^{-1} \|\bnabla \uInitPDE\|_2^{-4}$, hence relations \eqref{jwekrnwer}--\eqref{maxE} cannot rule out a singularity.

{\tolerance=1000%allows better line breaking
Numerical investigations of possible singularity formation in Navier-Stokes and Euler flows, the latter governed by system \eqref{nse}--\eqref{incompressible} with $\nu = 0$, based on a direct integration of these systems have had a long history and we refer the reader to the survey \cite{protas22} for a more in-depth discussion. However, in the context of Navier-Stokes flows, these efforts have been largely inconclusive. Here we follow a fundamentally different approach to searching for extreme, possibly singular, solutions which relies on variational optimization formulations where functionals related to different a priori estimates are maximized under suitable constraints. This research direction was initiated in \cite{luDoering08} where the authors considered the sharpness of bound \eqref{jwekrnwer}. An expression for the rate of growth of the enstrophy is obtained testing \eqref{nse} with $\Delta \u$, i.e. multiplying \eqref{nse} by $\Delta \u$, integrating the resulting expression over the domain $\Omega$ and then integrating by parts which gives
\begin{align}
    \G(\u) := \frac{d}{dt} \|\bnabla \u \|_2^2 = -\nu \int_{\Omega} | \Delta \u |^2 \, d\x + \int_{\Omega} \u \cdot \bnabla\u \cdot \Delta \u \, d\x.
    %\label{R}
\end{align}
}%
This then leads to the following 
\begin{problem}
    \label{pb:maxdEdt}
    Given $\E_0 > 0$, find 
    \begin{align}
        \tuE = \argmax_{\u\in \S_{\E_0}} \G(\u),
    \end{align}
    where
    \begin{align}
        \S_{\E_0} = \left\lbrace \u\in H^{2}(\Omega) \ \middle\vert\ \|\bnabla \u\|_2^2 = \E_0,\quad \bnabla \cdot \u = 0,\quad \int_{\Omega} \u \, d\x = \0\right\rbrace.
    \end{align}
\end{problem}
This problem was solved numerically for different $\E_0$ in \cite{luDoering08}, and later revisited in \cite{ayalaProtas2017}, which showed that $\G(\tuE) = C\E_0^3$ as $\E_0 \rightarrow \infty$ providing evidence for the sharpness of bound \eqref{jwekrnwer} (up to the numerical prefactor); this suggests the bound cannot be fundamentally improved. The optimal states $\tuE$ saturating this bound have the form of two colliding nearly-axisymmetric vortex rings which become increasingly localized as $\E_0 \rightarrow \infty$. The goal of the present study is to advance this research program and probe the sharpness of bound \eqref{estimate} for different $q > 3$, which will lead to variational optimization problems similar to Problem \ref{pb:maxdEdt}.

As regards the finite-time behavior, it was shown in \cite{ayalaProtas2017} that for very short times only the Navier-Stokes flows with the initial condition $\tuE$ saturate the inequality 
\begin{align}
    %\frac{1}{\|\bnabla \uInitPDE\|_2^2}-\frac{1}{\|\bnabla \u\|_2^2} \leq \frac{27}{4(2\pi\nu)^4} (\|\uInitPDE\|_2^2 - \|\u\|_2^2).\label{kwenrjewknr}
    \frac{1}{\|\bnabla \uInitPDE\|_2^2}-\frac{1}{\|\bnabla \u(t)\|_2^2} \leq \frac{27}{4(2\pi\nu)^4} (\|\uInitPDE\|_2^2 - \|\u(t)\|_2^2)\label{Et}
\end{align}
which is obtained by integrating \eqref{jwekrnwer} and using the energy equation
\begin{align}
    \frac{d}{dt}\|\u(t)\|_2^2 = -2 \nu \|\bnabla \u(t)\|_2^2,
    \label{dKdt}
\end{align}
which itself follows by testing \eqref{nse} with $\u$. On the other hand, for longer times, the nonlinear amplification is depleted leading to a very modest only growth of enstrophy. Thus, the main conclusion from \cite{luDoering08,ayalaProtas2017} was that if these are the only states saturating \eqref{jwekrnwer} and singularities do indeed form in Navier-Stokes flows, they are likely triggered by initial conditions that lead to a significant growth of the enstrophy but short of saturating bound \eqref{jwekrnwer}. Such states were then sought by maximizing the enstrophy, i.e. $\| \bnabla \u(t) \|_2^2$ at time $t>0$ and the integral in \eqref{lps} over a finite time with respect to the initial condition $\uInitPDE$ \cite{kangYunProtas2020,kangProtas2022,ramirez2025}. While no evidence was found for an unbounded growth of these quantities that would have signaled singularity formation, in all cases they exhibit a significant transient growth described by the relation 
\begin{align}
    \max_{t>0}\|\bnabla \u(t)\|_2^2 = C \|\bnabla \uInitPDE\|_2^3
    \label{E3/2}
\end{align}
in the limit $\|\bnabla \uInitPDE\|_2 \rightarrow \infty$.

We add that the methodology based on variational optimization formulations has also been employed to probe the sharpness of a priori bounds for some problems known to be globally well posed where these bounds are finite. The main interest here is that many of these estimates are obtained using analysis techniques similar to those used to obtain the bounds discussed above, so the question about their sharpness is also quite pertinent.  

In the context of one-dimensional (1D) viscous Burgers flows, the following bound on the rate of growth of 1D enstrophy was obtained in \cite{luDoering08}
\begin{align}
    \frac{d}{dt}\|\partial_x u(t)\|_2^2 \leq \|\partial_x u(t)\|_2^\frac{10}{3}
\end{align}
and was also shown to be achieved by states found by solving a 1D counterpart of Problem \ref{pb:maxdEdt}. As regards the maximum growth of enstrophy in finite time, the results from \cite{ayalaProtas2011} suggests that 
\begin{align}
    \max_{t>0}\| \partial_x u(t)\|_2^2 = C \|\partial_x u^0\|_2^3
\end{align}
in the limit $\|\partial_x u^0\|_2 \rightarrow \infty$ and a rigorous bound exhibiting this scaling was established only recently in \cite{AlbrittonDeNitti2023}.

Since in the absence of solid boundaries the enstrophy is a non-increasing quantity in two-dimensional (2D) Navier-Stokes flows, the relevant quantity in such situations is the palinstrophy $\| \Delta \u(t) \|_2^2$. The following a priori instantaneous bounds 
\begin{align}
    %\frac{d}{dt}\|\Delta \u\|_2^2 &\leq - \nu \frac{\|\Delta \u\|_2^4}{\|\bnabla \u\|_2^2} + \frac{C}{\nu}\|\bnabla \u\|_2^2 \|\Delta \u\|_2^2
    %\\
    \frac{d}{dt}\|\Delta \u(t)\|_2^2 &\leq - \nu \frac{\|\Delta \u(t)\|_2^4}{\|\bnabla \u(t)\|_2^2} + \frac{C}{\nu}\|\bnabla \u(t)\|_2^2 \|\Delta \u(t)\|_2^2,
    \label{ayalaProtas2014_1}
    \\
    %\frac{d}{dt}\|\Delta \u\|_2^2 &\leq \frac{C}{\nu} \|\bnabla \u\|_2 \|\Delta \u\|_2^3
    %\\
    \frac{d}{dt}\|\Delta \u(t)\|_2^2 &\leq \frac{C}{\nu} \|\bnabla \u(t)\|_2 \|\Delta \u(t)\|_2^3,
    \label{ayalaProtas2014_2}
\end{align}
and finite-time estimates
\begin{align}
    %\max_{t>0} \|\Delta \u\|_2^2 &\leq \|\Delta \uInitPDE\|_2^2 + \frac{C}{\nu^2}\|\bnabla \uInitPDE\|_2^4
        \max_{t>0} \|\Delta \u(t)\|_2^2 &\leq \|\Delta \uInitPDE\|_2^2 + \frac{C}{\nu^2}\|\bnabla \uInitPDE\|_2^4,
    \\
    %\max_{t>0} \|\Delta \u\|_2^2 &\leq \left(\|\Delta \uInitPDE\|_2 + \frac{C}{\nu^2}\|\bnabla \uInitPDE\|_2\|\Delta \uInitPDE\|_2^2\right)^2
    \max_{t>0} \|\Delta \u(t)\|_2^2 &\leq \left(\|\Delta \uInitPDE\|_2 + \frac{C}{\nu^2}\|\bnabla \uInitPDE\|_2\|\Delta \uInitPDE\|_2^2\right)^2
    \label{ayalaProtas2014_4}
\end{align}
on, respectively, the rate of growth and the finite-time behavior of this quantity were established in \cite{ayalaProtas2014}, where they were also shown to be sharp with respect to different parameters. Bound \eqref{ayalaProtas2014_2} was revisited in \cite{ayalaDoeringSimon2018} where an estimate with an improved prefactor was established
\begin{align}
    %\frac{d}{dt}\|\Delta \u\|_2^2 &\leq C\sqrt{\ln \frac{\|\bnabla \u \|_2}{\nu}+C}\ \|\Delta \u\|_2^3.
    \frac{d}{dt}\|\Delta \u(t)\|_2^2 &\leq C\sqrt{\ln \frac{\|\bnabla \u(t) \|_2}{\nu}+C}\ \|\Delta \u(t)\|_2^3.
\end{align}

A condensed summary of these results is given in Table \ref{table_previous_results} with additional details provided in the survey paper \cite{protas22}. The present study complements this research program by providing computational evidence that bound \eqref{estimate} is in fact sharp. In addition, it is also argued that the states saturating these bounds are far from saturating estimate \eqref{jwekrnwer} on the rate of growth of enstrophy.

\begin{table}[ht]
    \centering
    \resizebox{\linewidth}{!}{
\newcommand{\myGap}{\\\\}
\newcommand{\myHspace}{\hspace{20pt}}
\newcommand{\myMultiRow}[1]{\multirow{2}{*}{#1}}
\newcommand{\myHeaderGap}{-2pt}
\newcommand{\myCategoryGap}{10pt}
\newcommand{\myTabSep}{10pt}
$%
    \begin{alignedat}{3}
        &\text{Problem} & \clap{\text{Estimate}}& & &\text{Numerically Realized}
        %\\[-15pt]
        \\
        \hline
        %\intertext{\hrulefill}
        \\[-10pt]
        &\text{Navier-Stokes equations}
        \hspace{\myTabSep} & \hspace{\myTabSep}
        \\[\myHeaderGap]
        &\myHspace\text{2d, Instantaneous}
        \hspace{\myTabSep} & \hspace{\myTabSep} 
        \frac{d}{dt}\|\Delta \u\|_2^2 &\lesssim \|\Delta \u\|_2^3
        \hspace{\myTabSep} & \hspace{\myTabSep}
        &\text{Yes \cite{ayalaProtas2014,ayalaDoeringSimon2018}}
        \\
        &\myHspace\text{2d, Finite Time}
        \hspace{\myTabSep} & \hspace{\myTabSep} 
        \max\limits_{t>0}\|\Delta \u\|_2^2 &\lesssim \|\Delta \u_0\|_2^2+\|\bnabla \u_0\|_2^4
        \hspace{\myTabSep} & \hspace{\myTabSep} 
        &\text{Yes \cite{ayalaProtas2014,ayalaDoeringSimon2018}}
        \\
        &\myHspace\text{3d, Instantaneous}
        \hspace{\myTabSep} & \hspace{\myTabSep} 
        \frac{d}{dt}\|\bnabla \u\|_2^2 &\lesssim \|\bnabla \u\|_2^6
        \hspace{\myTabSep} & \hspace{\myTabSep}
        &\text{Yes \cite{ayalaProtas2017,luDoering08}}
        \\
        &\myHspace\text{}
        \hspace{\myTabSep} & \hspace{\myTabSep} 
        \frac{d}{dt}\|\u\|_q^q &\lesssim \|\u\|_q^{q\frac{q-1}{q-3}}
        \hspace{\myTabSep} & \hspace{\myTabSep}
        &\text{Yes [present work]}
        \\
        &\myHspace\text{3d, Finite Time}
        \hspace{\myTabSep} & \hspace{\myTabSep} 
        \|\bnabla \u\|_2^2 &\leq \frac{\|\bnabla \u_0\|_2^2}{\sqrt{1-c\|\bnabla \u_0\|_2^4t}}
        \hspace{\myTabSep} & \hspace{\myTabSep}
        &\text{No \cite{ayalaProtas2017,kangYunProtas2020}}
        \\
        &\myHspace\text{}
        \hspace{\myTabSep} & \hspace{\myTabSep} 
        \int_0^T \|\u\|_4^\frac{8}{3} \, dt &\lesssim \|\u_0\|_2^\frac{8}{3}
        \hspace{\myTabSep} & \hspace{\myTabSep}
        &\text{No \cite{kangProtas2022}}
        \\
        &\myHspace\text{}
        \hspace{\myTabSep} & \hspace{\myTabSep} 
        \int_0^T \|\u\|_q^{\frac{2q}{q-3}}\, dt &\to \infty
        \hspace{\myTabSep} & \hspace{\myTabSep}
        &\text{No \cite{ramirez2025,kangProtas2022}}
        \\[\myCategoryGap]
        &\text{Euler equations}
        \hspace{\myTabSep} & \hspace{\myTabSep}
        \\[\myHeaderGap]
        &\myHspace\text{3d, Finite Time}
        \hspace{\myTabSep} & \hspace{\myTabSep} 
        \max\limits_{0\leq t\leq T}\|\u\|_{H^3} &\to \infty \phantom{\frac{d}{dt}}
        \hspace{\myTabSep} & \hspace{\myTabSep}
        &\text{Yes \cite{zhaoProtas2023}}
        \\[\myCategoryGap]
        &\text{Viscous Burgers' equations}
        \hspace{\myTabSep} & \hspace{\myTabSep}
        \\[\myHeaderGap]
        &\myHspace\text{1d, Instantaneous}
        \hspace{\myTabSep} & \hspace{\myTabSep} 
        \frac{d}{dt}\|\partial_x u\|_2^2 &\lesssim \|\partial_x u_0\|_2^\frac{10}{3}
        \hspace{\myTabSep} & \hspace{\myTabSep}
        &\text{Yes \cite{luDoering08}}
        \\
        &\myHspace\text{1d, Finite Time}
        \hspace{\myTabSep} & \hspace{\myTabSep} 
        \max\limits_{0\leq t\leq T} \|\partial_x u\|_2^2 &\lesssim \|\partial_x u_0\|_2^3
        \hspace{\myTabSep} & \hspace{\myTabSep}
        &\text{Yes \cite{ayalaProtas2011}}
    \end{alignedat}
$%
}
    \caption{An overview of selected a priori bounds and regularity conditions for hydrodynamic models in 1D, 2D and 3D together with information about their realizability in numerical computations.}%
    \label{table_previous_results}%
\end{table}

%\subsubsection{outline of the paper}

The structure of the paper is as follows. In Section \ref{sectionOptimizationProblemDerivation} we introduce the objective functional as well as the maximization problems corresponding to \eqref{estimate} in the setting of Banach and Hilbert spaces. In Section \ref{sectionSmallDataLimit} we analytically solve the problem in the small-data limit producing solutions used as the initial guesses for the numerical approach described in Section \ref{sectionNumericalApproach}. The results are presented in Section \ref{sectionResults} and conclusions are deferred to Section \ref{sectionConclusion}.

%\subsubsection{notation}
\paragraph{Notation}\leavevmode\\
Throughout the paper, $L^q$ denotes the Lebesgue space of $q$-integrable functions equipped with the norm $\|\cdot\|_q$ defined by
\begin{alignat}{2}
    &\|\boldsymbol{f}\|_q^q && \defeq \int_\Omega |\boldsymbol{f}(\x)|^q\ d\x,
    \qquad \text{for} \ 1\leq q <\infty, \\
    &\|\boldsymbol{f}\|_\infty && \defeq \esssup_{\x\in \Omega} |\boldsymbol{f}(\x)|,
\end{alignat}
$W^{k,q}$ is the Sobolev space of functions with $q$-integrable weak derivatives of order up to $k$ equipped with the norm $\|\cdot\|_{W^{k,q}}$ defined by
\begin{align}
    %{\color{gray}\|\boldsymbol{f}\|_{k,q}^q \defeq \sum_{|\boldsymbol{\alpha}|\leq k}\int_\Omega |\bnabla^{\boldsymbol{\alpha}} \boldsymbol{f}(\x)|^q\ d\x}
    %\\
    \|\boldsymbol{f}\|_{W^{k,q}}^q \defeq \int_\Omega \left|\left(1+\bnabla^{k}\right) \boldsymbol{f}(\x)\right|^q\ d\x
\end{align}
and $H^k\defeq W^{k,2}$ is a Hilbert space equipped with the same norm. 
Similarly, we define the corresponding inner products
\begin{align}
    \langle \boldsymbol{f},\boldsymbol{g}\rangle_{L^2} &\defeq \int_{\Omega} \boldsymbol{f}(\x)\cdot\boldsymbol{g}(\x)\, d\x,
    \\
    %&{\color{gray}\langle \boldsymbol{f},\boldsymbol{g}\rangle_{H^k} \defeq \sum_{|\boldsymbol{\alpha}|\leq k} \int \bnabla^{\boldsymbol{\alpha}} \boldsymbol{f}(\x)\mathcolon \bnabla^{\boldsymbol{\alpha}}\boldsymbol{g}(\x)\, d\x}\\
    \langle \boldsymbol{f},\boldsymbol{g}\rangle_{H^s} &\defeq \sum_{j=0,s} \int_{\Omega} \bnabla^{j} \boldsymbol{f}(\x)\mathcolon\bnabla^{j} \boldsymbol{g}(\x)\, d\x. \label{Hs}
\end{align}

\section{Optimization Problems}
\label{sectionOptimizationProblemDerivation}

The goal of this research is to probe the sharpness of bound \eqref{estimate} in the limit $\| \u \|_q \rightarrow \infty$. To do so, we will formulate an optimization where an expression for $(d/dt)\| \u(t) \|_q$ is maximized with respect to the vector field $\u$ subject to constraints which include fixing the norm $\| \u \|_q = B$. The objective is to check whether the maximum rate of growth of the $L^q$ norm determined by solving this optimization problem saturates the upper bound in estimate \eqref{estimate} as $B \rightarrow  \infty$ and, if so, what vector fields $\tuB$ realize this behavior. Formulation of this optimization problem requires specification of the objective functional $\R_q(\u)$ together with an appropriate solution space and a constraint, which is done in the subsections below.

\subsection{Objective Functional}
\label{subsectionObjectiveFunctional}

To obtain an expression for $(d/dt)\| \u(t) \|_q$ we begin by testing \eqref{nse} with $|\u|^{q-2}\u$ which yields
\begin{align}
    \frac{1}{q}\frac{d}{dt}\|\u\|_q^q &= \int_{\Omega} |\u|^{q-2} \u\cdot \partial_t \u \, d\x
    \\
    &= - \int_{\Omega} |\u|^{q-2} \u\cdot (\u\cdot\bnabla) \u \, d\x - \int_{\Omega} |\u|^{q-2} \u\cdot \bnabla p \, d\x
    \nonumber \\
    &\qquad + \nu \int_{\Omega} |\u|^{q-2}\u\cdot \Delta \u \, d\x.
    \label{dLqdt}
\end{align}
Note that the first term on the right-hand side of \eqref{dLqdt} vanishes due to incompressibility condition \eqref{incompressible} and the periodic boundary conditions under integration by parts as
\begin{align*}
    0 &= \int_{\Omega} |\u|^{q} \bnabla \cdot \u \, d\x = - \int_{\Omega} \u\cdot \bnabla \left(|\u|^{q}\right) \, d\x = -q \int_{\Omega} |\u|^{q-2} \u\cdot (\u\cdot\bnabla) \u \, d\x.
\end{align*}
Expression \eqref{dLqdt} depends on both velocity $\u$ and pressure $p$, however, it will be more convenient if the functional in the optimization problem depends on $\u$ alone. Taking the divergence of \eqref{nse} and using \eqref{incompressible} gives the Poisson equation satisfied by pressure
\begin{align}
    \Delta p = - \bnabla \u\mathcolon\bnabla \u^T, \qquad \int_{\Omega} p \, d\x = 0, \label{p}
\end{align}
where the zero-mean condition fixes the arbitrary additive constant. The term involving pressure on the right-hand side of \eqref{dLqdt} can be therefore transformed using integration by parts and equation \eqref{p} as
\begin{align}
    \int_\Omega |\u|^{q-2} \u\cdot \bnabla p\, d\x &= - \int_\Omega |\u|^{q-2} \bnabla \cdot \u\, p \, d\x - \int_\Omega \u\cdot \bnabla |\u|^{q-2}p\, d\x
    \\
    &= - (q-2)\int_\Omega |\u|^{q-4} \u\cdot (\u\cdot \bnabla) \u \, p \, d\x
    \\
    &= - (q-2)\int_\Omega |\u|^{q-4} \u\cdot (\u\cdot \bnabla) \u\, \Delta^{-1}(\bnabla \u\mathcolon \bnabla \u^T)\, d\x.\quad\label{omittedPressure}
\end{align}
After performing integration by parts the last term on the right-hand side of \eqref{dLqdt} becomes
\begin{align}
     \int_{\Omega} |\u|^{q-2}\u\cdot \Delta \u \, d\x &= - \int_{\Omega} |\u|^{q-2}|\bnabla \u|^2 \, d\x
     \nonumber \\
     &\qquad - (q-2)\int_{\Omega} |\u|^{q-4}u_i \partial_j u_i u_k \partial_j u_k\, d\x 
     \\
     &= - \int_{\Omega} |\u|^{q-2}|\bnabla \u|^2 \, d\x - \frac{4(q-2)}{q^2}\int_{\Omega} \left|\bnabla |\u|^{\frac{q}{2}}\right|^2 d\x.
     \quad\label{dLqdt3}
\end{align}
Therefore, combining \eqref{dLqdt}, \eqref{omittedPressure} and \eqref{dLqdt3}, we arrive at the following expression for the objective functional that depends on $\u$ only
\begin{align}
    \frac{1}{q}\frac{d}{dt}\|\u\|_q^q &= \int_\Omega |\u|^{q-2} \u\cdot \bnabla \Delta^{-1} (\bnabla \u\mathcolon\bnabla \u^T) \, d\x + \nu \int_{\Omega} |\u|^{q-2}\u\cdot \Delta \u \, d\x
    \label{sakfndsj}
    \\
    &= - (q-2)\int_\Omega |\u|^{q-4} \u\cdot (\u\cdot \bnabla) \u\, \Delta^{-1}(\bnabla \u\mathcolon \bnabla \u^T)\, d\x \nonumber
    \\
    &\qquad - \nu  \int_{\Omega} |\u|^{q-2}|\bnabla \u|^2 \, d\x - \frac{4(q-2)\nu}{q^2}\int_{\Omega} \left|\bnabla |\u|^{\frac{q}{2}}\right|^2d\x \nonumber
    \\
    &\eqdef \R_q(\u).
    \label{defRq}
\end{align}

\subsection{Functional Setting}
\label{subsectionFunctionalSetting}

Here we determine the minimum regularity of the argument of $\R_q(\u)$ required for this functional to be well defined. Using H{\"o}lder's inequality the objective functional in the form \eqref{dLqdt} can be bounded as
\begin{align}
    |\R_{q}(\u)| &\leq C \nu \|\u\|_{3q}^{q-2} \|\bnabla \u\|_{\frac{3q}{q+1}}^2 + C\|p\|_{\frac{3(q+1)}{4}}\|\u\|_{\frac{3(q+1)}{2}}^{q-2} \|\bnabla \u\|_{\frac{3(q+1)}{q+3}},
    \label{Rqbound}
\end{align}
where the constant $C>0$ only depends on $q$. The pressure, as a solution of \eqref{p}, satisfies the bound (for details, see Lemma 3 in \cite{robinsonSadowski2014})
\begin{align}
    \|p\|_{s}\leq \|\u\|_{2s}^2 \qquad \text{for all} \ s > 1. \label{pu2}
\end{align}
Applying \eqref{pu2} and a Sobolev embedding \cite{af05} to \eqref{Rqbound} yields
\begin{align}
    |\R_{q}(\u)| &\leq C \nu \|\u\|_{3q}^{q-2} \|\bnabla \u\|_{\frac{3q}{q+1}}^2 + C\|\u\|_{\frac{3(q+1)}{2}}^{q} \|\bnabla \u\|_{\frac{3(q+1)}{q+3}}
    \\
    &\leq C \nu \|\u\|_{W^{1,\frac{3q}{q+1}}}^q + C \|\u\|_{W^{1,\frac{3(q+1)}{q+3}}}^{q+1}
    \\
    &\leq C \|\u\|_{W^{1,\frac{3q}{q+1}}}^q \left(\nu + \|\u\|_{W^{1,\frac{3q}{q+1}}}\right),
    \label{bjwerbewhr}
\end{align}
where in the last step we again used H{\"o}lder's inequality since
\begin{align}
    \frac{3(q+1)}{q+3} \leq \frac{3q}{q+1}
\end{align}
for $q\geq 1$ and we consider $q\geq 3$. Therefore, in order to have a well-defined objective functional, we require the velocity field to have the regularity $\u\in W^{1,\frac{3q}{q+1}}(\Omega)$. For reasons that will become clear in Section \ref{sectionNumericalApproach}, we also need to consider a formulation in a Hilbert space, which we choose to be the largest Sobolev-Hilbert space $H^s(\Omega)$ embedded in $W^{1,\frac{3q}{q+1}}(\Omega)$, as determined by the Sobolev embedding
\begin{align}
    H^{\frac{3}{2}-\frac{1}{q}}(\Omega) \hookrightarrow W^{1,\frac{3q}{q+1}}(\Omega), \qquad q \ge 3.
    \label{HsWp}
\end{align}

\subsection{Statement of the Optimization Problems}

The considerations in Sections \ref{subsectionObjectiveFunctional} and \ref{subsectionFunctionalSetting} immediately motivate the following two variational problems.
\begin{problem}[Banach-space Problem]
    \label{banachProblem}
    Let $\R_q$ be defined by \eqref{defRq} with $q>3$ and $B>0$ be given, find
    \begin{align}
        \tuB = \argmax_{\u\in Y_B} \R_q(\u),
    \end{align}
    where
    \begin{align}
        Y_B = \left\lbrace \u\in W^{1,\frac{3q}{q+1}}\ \middle\vert\ \|\u\|_q=B,\quad \bnabla \cdot \u = 0,\quad \int_{\Omega} \u \, d\x = \0\right\rbrace.
    \end{align}
\end{problem}
\begin{problem}[Hilbert-space Problem]
    \label{hilbertProblem}
    Let $\R_q$ be defined by \eqref{defRq} with $q>3$ and $B>0$ be given, find
    \begin{align}
        \tuB = \argmax_{u\in X_B} \R_q(\u),
    \end{align}
    where
    \begin{align}
        X_B = \left\lbrace \u\in H^{\frac{3}{2}-\frac{1}{q}} \ \middle\vert\ \|\u\|_q=B,\quad \bnabla \cdot \u = 0,\quad \int_{\Omega} \u \, d\x = \0\right\rbrace.
    \end{align}
\end{problem}
The norm constraint $\|\u\|_q$ is needed as in its absence the objective functional could be made arbitrarily large simply by increasing the magnitude of the argument $\u$, whereas the zero-mean condition excludes solutions in the form of constant vector fields. The functional setting of Problem \ref{hilbertProblem} naturally endows the manifold $X_B$ with an inner product thus making the problem Riemannian. This structure will be exploited in the design of our solution approach described in Section \ref{sectionNumericalApproach}.

Since the space \H{} is embedded in $W^{1,\frac{3q}{q+1}}$, cf.~\eqref{HsWp},  if a maximizing sequence of Problem \ref{hilbertProblem} shows that bound \eqref{estimate} is sharp, then solving Problem \ref{banachProblem} becomes unnecessary and can be omitted. As Problem \ref{hilbertProblem} is both mathematically and computationally more tractable and our results obtained by solving this problem do in fact suggest that bound \eqref{estimate} is sharp, hereafter we primarily focus on this case.

%\subsection{Existence of Maximizers}
%\todo[inline]{The proof of the existence needs rewriting and some ideas/discussion. Uncomment input\{existenceOfMaximizers\} to see the unfinished proof.}
%\input{existenceOfMaximizers}

\section{Solution of Problems \ref{banachProblem} and \ref{hilbertProblem} in the Small-Data Limit}
\label{sectionSmallDataLimit}
%\subsubsection{small-data limit \texorpdfstring{$\to$}{to} eigenvalue problem}

In this section we construct (some of the) solutions of Problems \ref{banachProblem} and \ref{hilbertProblem} in the small-data limit following the approach developed in \cite{ayalaProtas2014,ayalaProtas2017}. 
The construction works for both, Problem \ref{banachProblem} and \ref{hilbertProblem}. In addition to their independent interest, these maximizers will be used as initial guesses to facilitate numerical solution of Problem \ref{hilbertProblem} for finite values of the parameter $B$, cf.~Section \ref{sectionResults}. Intuitively, when $B \rightarrow 0$, the pressure term in the objective functional \eqref{dLqdt} is dominated by the viscous term, such that maximizers in this limit include eigenfunctions of the Laplacian. To show this, we convert Problem \ref{banachProblem} to an unconstrained form with the Lagrangian 
\begin{align}
    \mathcal{L}(\u,\lambda,r) &\defeq \mathcal{R}_{q}(\u) + \int_\Omega r \bnabla \cdot \u\ d\x + \lambda (\|\u\|_q^q-B^q)
    \\
    &=- \nu \left\| |\u|^{\frac{q-2}{2}}|\bnabla \u|\right\|_2^2 - \frac{4(q-2)}{q^2} \nu \left\|\bnabla|\u|^{\frac{q}{2}}\right\|_2^2 \nonumber
    \\
    &\qquad + (q-2) \int_\Omega p |\u|^{q-4} \u\cdot (\u\cdot \bnabla) \u + \int_\Omega r \bnabla \cdot \u\ d\x \nonumber
    \\
    &\qquad + \lambda (\|\u\|_q^q - B^q),
\end{align}
where $\lambda \in \RR$ and $r \in L^\frac{3q}{2q-1}(\Omega)$ are Lagrange multipliers.
The corresponding Euler-Lagrange equations obtained by differentiating the Lagrangian with respect to its three arguments are
\begin{align}
    (q-2)|\u|^{q-4} \u \u\cdot (\nu\Delta \u - \bnabla p) + |\u|^{q-2}(\nu\Delta \u- \bnabla p)\quad & \nonumber
    \\
    - 2 \bnabla (\u\cdot \bnabla) \Delta^{-1} \bnabla \cdot (|\u|^{q-2}\u) + \nu \Delta(|\u|^{q-2}\u) - \bnabla r + \lambda q |\u|^{q-2}\u &= \0,\quad\label{EL1_Small}
    \\
    \bnabla \cdot \u &= 0,\quad\label{EL2_Small}
    \\
    \|\u\|_q - B &= 0.\quad\label{EL3_Small}
\end{align}
We then make the following ans{\"a}tze for some $\alpha>0$ which is to be determined
\begin{alignat}{4}
    \u &= \u_0 &&+ B^\alpha \u_1 &&+ B^{2\alpha}\u_2 &&+ \dots, \label{au}\\
    \lambda &= \lambda_0 &&+ B^\alpha \lambda_1 &&+ B^{2\alpha}\lambda_2 &&+ \dots, \label{al} \\
    r &= r_0 &&+ B^\alpha r_1 &&+ B^{2\alpha}r_2 &&+ \dots, \label{ar}
\end{alignat}
and insert them in \eqref{EL1_Small}--\noeqref{al}\eqref{EL3_Small}. The last of these equations then becomes 
\begin{align}
    B^q &= \|\u\|_q^q
    \\
    &= \|\u_0\|_q^q + B^\alpha \int_\omega |\u_0|^{q-2}\u_0\cdot \u_1\, d\x + B^{2\alpha}\int_\Omega |\u_0| (1+\mathcal{O}(\u_0,B^\alpha))\, d\x \nonumber
    \\
    &\qquad + B^{\alpha q}\|\u_1\|_q^q + B^{\alpha(q+1)}q\int_\Omega |\u_1|^{q-2}\u_1\cdot\u_2 \, d\x+ \mathcal{O}(B^{\alpha(q+2)}).
    \label{lkjasdf}
\end{align}
As $B \rightarrow 0$, at the leading (zeroth) order we find
\begin{align}
    \|\u_0\|_q^q = 0
\end{align}
which implies that the first non-trivial correction has to satisfy
\begin{align}
    B^q = B^{\alpha q}\|\u_1\|_q^q,
\end{align}
and therefore $\alpha = 1$ and $\|\u_1\|_q = 1$. At the lower orders we find $\bnabla r_k=0$, implying that $r_k=0$ for all $k<q-1$.  The lowest nontrivial correction in \eqref{EL1_Small} is given at the order $B^{q-1}$ by
\begin{align}
    (q-2)\nu|\u_1|^{q-4} \u_1 \u_1\cdot \Delta \u_1 + \nu|\u|^{q-2}\Delta \u_1\qquad\qquad\qquad & \nonumber
    \\
    + \nu \Delta(|\u_1|^{q-2}\u_1) - \bnabla r_{q-1} + \lambda_0 q |\u_1|^{q-2}\u_1 &= \0.\label{u1}
\end{align}
Testing this relation with $\u_1$ we find
\begin{align}
    0 &= (q-1)\nu\int_\Omega |\u_1|^{q-2} \u_1 \cdot \Delta \u_1 \, d\x + \nu\int_\Omega \u_1 \cdot \Delta (|\u_1|^{q-2}\u_1) d\x \nonumber
    \\
    &\qquad - \int_\Omega \u_1\cdot \bnabla r_{q-1}\, d\x + \lambda_0 q \|\u_1\|_q^q
    \\
    &= q \int_\Omega |\u_1|^{q-2}\u_1\cdot (\nu \Delta \u_1 + \lambda_0 \u_1) d\x,
    \label{lkwerjkwej}
\end{align}
where in the second equality we used integration by parts and $\bnabla\cdot\u_k=0$ due to \eqref{EL2_Small}. Evaluating the objective functional using the ansätze \eqref{au}--\eqref{ar} we find
\begin{align}
    \frac{1}{q}\frac{d}{dt}\|\u\|_q^q &= \R_q(\u) =  \int_\Omega |\u|^{q-2} \u\cdot \bnabla p \, d\x + \nu \int_{\Omega} |\u|^{q-2}\u\cdot \Delta \u \, d\x
    \\
    &= -\lambda_0 B^q \|\u_1\|_q^q + \mathcal{O}(B^{q+1})= -\lambda_0 B^q + \mathcal{O}(B^{q+1}).\label{objectiveExpansion1}
\end{align}
Trivially, eigenfunctions of the Laplacian satisfy \eqref{lkwerjkwej}, provided the Lagrange multiplier $\lambda_0$ satisfies $\Lambda = \lambda_0 / \nu$, where $\Lambda$ is the corresponding eigenvalue. Here we use the standard convention $-\Delta u = \Lambda u$ for eigenfunctions $u$ and eigenvalues $\Lambda$ of the Laplacian. However, any sufficiently smooth, divergence-free, nontrivial $\u_1$ is also a solution of \eqref{lkwerjkwej} for some $\lambda_0\in \RR$. In fact, we note that $\|u_1\|_q^q>0$ and by \eqref{dLqdt3}
\begin{align*}
    \nu\int_\Omega |\u_1|^{q-2}\u_1\cdot \Delta \u_1\, d\x = - \nu \left\| |\u_1|^{\frac{q-2}{2}}|\bnabla \u_1|\right\|_2^2 - \frac{4(q-2)}{q^2} \nu \left\|\bnabla|\u_1|^{\frac{q}{2}}\right\|_2^2 <0,
\end{align*}
such that \eqref{lkwerjkwej} is satisfied for
\begin{align}
    \lambda_0 = \frac{ -\nu \int_\Omega |\u_1|^{q-2}\u_1\cdot \Delta \u_1\, d\x}{\|\u_1\|_q^q}.
\end{align}
Expression \eqref{objectiveExpansion1} is the maximized for the modified Rayleigh quotient
\begin{align}
    \widetilde \lambda_0 = \min_{\u_1} \frac{-\nu \int_\Omega |\u_1|^{q-2}\u_1\cdot \Delta \u_1\, d\x }{\|\u_1\|_q^q}.
\end{align}
However, since the eigenfunctions $\v_k$ of the Laplacian form an orthogonal basis of $L^2(\Omega)$ and can be ordered by the corresponding eigenvalues $0<\Lambda_1 < \Lambda_2 <\dots$, we find
\begin{align}
    -\int_\Omega |\u_1|^{q-2}\u_1\cdot \Delta \u_1\, d\x &= -\sum_{k=1}^\infty \int_\Omega |\u_1|^{q-2} \u_1\cdot c_k\Delta \v_k\, d\x
    \\
    &= \sum_{k=1}^\infty \int_\Omega |\u_1|^{q-2} \u_1\cdot \Lambda_k c_k\v_k\, d\x
    \\
    &\geq \min_k \Lambda_k \|\u_1\|_q^q
    \\
    &= \Lambda_1 \|\u_1\|_q^q,
\end{align}
with equality if and only if $\u_1$ is the eigenfunction of the Laplacian corresponding to the smallest eigenvalue $\Lambda_1$, implying $\widetilde \lambda_0 = \nu \Lambda_1$.

Therefore, in the limit $B \rightarrow 0$, the objective functional $\R_q(\u)$ is maximized by the eigenfunction of the Laplacian corresponding to the smallest eigenvalue $\Lambda_1$ and, in view of \eqref{objectiveExpansion1}, it is given by
\begin{align}
    \frac{1}{q}\frac{d}{dt}\|\u\|_q^q &= \R_q(\u) = - \nu \Lambda_1 B^q + \mathcal{O}(B^{q+1}).\label{objectiveExpansion2}
\end{align}
The eigenfunction corresponding to the smallest eigenvalue $\Lambda_1=4\pi^2$ is, up to shifts, given by Arnold-Beltrami-Childress flow \cite{ayalaProtas2017, mb02}
\begin{align}
    \u_1(\x) = \begin{pmatrix}
        A\sin(2\pi x_3) + C\cos(2\pi x_2)\\
        B\sin(2\pi x_1) + A\cos(2\pi x_3)\\
        C\sin(2\pi x_2) + B\cos(2\pi x_1)
    \end{pmatrix},
    \label{ABC}
\end{align}
which yields
\begin{align}
    \frac{1}{q}\frac{d}{dt}\|\u\|_q^q &= \R_q(\u) = - 4\pi^2 \nu  B^q + \mathcal{O}(B^{q+1}).
\end{align}
It will be used with $A=B=C=1$ as the initial guess $\uInit$ for the continuation approach described in Algorithm \ref{alg1} below, i.e., to initialize the computation of branches of local maximizers by successively solving Problem \ref{hilbertProblem} for increasing values of $B$.

%and serves as the initiation of the branches we calculate. Specifically, we will set $A=B=C=1$ and use
% \begin{align}
%     \u^{\text{init}}(\x) = \begin{pmatrix}
%         \cos(2\pi x_2) + \cos(2\pi x_3)\\
%         \cos(2\pi x_3) + \cos(2\pi x_1)\\
%         \cos(2\pi x_1) + \cos(2\pi x_2)
%     \end{pmatrix}
% \end{align}
%as our the initial guess.

\begin{remark}
    We note that \eqref{objectiveExpansion2} is strictly negative as in the limit considered the rate of growth of $\|\u(t)\|_q$ is dominated by the viscous (dissipative) effects. As will be confirmed by our results in Section \ref{sectionResults}, the functional \eqref{defRq} becomes positive for sufficiently large values of the constraint parameter $B$.
\end{remark}

\begin{remark}
    We note that constant solutions maximize the objective functional \eqref{defRq} for small values of $B$. The zero-mean constraint in Problems \ref{banachProblem} and \ref{hilbertProblem} excludes such trivial solutions.
\end{remark}

\section{Numerical Approach}
\label{sectionNumericalApproach}
%\subsubsection{Discrete Gradient Flow}
Here we describe our approach to numerically solve Problem \ref{hilbertProblem} which is an adaptation of the techniques developed earlier to solve optimization problems with a similar structure \cite{ayalaProtas2017,kangProtas2022,zhaoProtas2023}. For clarity, we first offer a high-level description of the approach and then fill in all the technical details. The reader is referred to the monograph \cite{absilMahonySepulchre} for additional information.

\subsection{Riemannian Conjugate-Gradient Approach}
\label{subsectionQualitativeDescription}%
Fundamentally, the strategy is based on a gradient flow for functional \eqref{defRq}
\begin{align}
    \frac{d\u}{d\tau} = \nabla\R_q(\u(\tau))
\end{align}
which after an explicit discretization with respect to the pseudo-time $\tau$ yields the iterative scheme
\begin{align}
    \u_{k+1} &=\u_k + \tau_k \nabla\R_q(\u_k), \qquad k = 0,1,2,\dots
\end{align}
for some step sizes $\tau_k>0$. Given some $\u_k \in X_B$, the new iterate $\u_{k+1}$ will in general not satisfy this constraint; thus, a retraction must be performed to ensure $\u_{k+1}\in X_B$, i.e.,
\begin{align}
    \u_{k+1} &=\Retraction\left(\u_k + \tau_k \nabla\R_q(\u_k)\right).
\end{align}
Given a suitable initial guess $\u_0 \in X_B$, it is expected that local maximizers of Problem \ref{hilbertProblem} can be approximated as $\tuB = \lim_{k \rightarrow \infty} \u_k$.

In order to accelerate the convergence, the gradient is projected onto the space $\TukXB$ tangent to the constraint manifold $X_B$ at $\u_k$, yielding the Riemannian gradient method
\begin{align}
    \u_{k+1} &=\Retraction\left(\u_k + \tau_k \mathbb{P}_k \nabla\R_q(\u_k)\right).
    \label{numExplRiemannScheme}
\end{align}
To further accelerate the procedure we employ the conjugate-gradient method which uses the information from the previous steps \cite{nw00}. In the unconstrained Euclidean setting it is given by
\begin{align}
    \u_{k+1} &=  \u_k + \tau_k  \d_k
    \label{numExplConjScheme_a}
    \\
    \d_k &=\nabla\R_q(\u_k) + \beta_k \d_{k-1}
    \label{numExplConjScheme_b}
\end{align}
where $\d_k$ is the step direction and $\beta_k \in \mathbb{R}$ is referred to as the momentum coefficient.

To combine both schemes, \eqref{numExplRiemannScheme} and \eqref{numExplConjScheme_a}--\eqref{numExplConjScheme_b}, the previous direction $\d_{k-1}$ which belongs to $\TukmoXB$, the tangent space at $\u_{k-1}$, has to be mapped to the new tangent space $\TukXB$ at $\u_k$ which is done via the vector transport $\Gamma$. This yields the Riemannian conjugate-gradient method
\begin{align}
    \u_{k+1} &=\Retraction\left(\u_k + \tau_k \d_k\right)
    \label{numExplRiemConjScheme_a}
    \\
    \d_k &=\nabla\R_q(\u_k) + \beta_k \Gamma(\d_{k-1}).
    \label{numExplRiemConjScheme_b}
\end{align}
The complete form of \eqref{numExplRiemConjScheme_a}--\eqref{numExplRiemConjScheme_b} is derived in the following subsections and given in \eqref{RCG1}--\eqref{RCG2}.

\subsection{Hilbert-Sobolev Gradient}
\label{subsectionHilberSobolevGradient}%
We note that the gradient $\nablaH\R_q(\u)$ needs to be defined as an element of the space \H. The starting point is the computation of the Gateaux derivative of $\R_q$ at $\v \in \text{\H{}}$ in some direction $\v' \in \text{\H{}}$ which is defined as $d\R_q(\v;\v') := \lim_{\epsilon \rightarrow 0} \epsilon^{-1}\left[ \R_q(\v + \epsilon \v') - \R_q(\v) \right]$. The calculations detailed in \ref{appendixCalcRderivative} give
\begin{align}
    d\R_q(\v;\v') &=  \int_\Omega |\v|^{q-4} \big[(q-2) (\v\cdot \v') \v + |\v|^2 \v' \big] \nonumber
    \\
    &\qquad \qquad \qquad \qquad \cdot \big[ \nu \Delta \v + \bnabla \Delta^{-1} \left(\bnabla \v\mathcolon\bnabla \v^T\right)\big] d\x \nonumber
    \\
    &\qquad + 2 \int_\Omega (\v' \cdot \bnabla) \v \cdot \bnabla  \Delta^{-1} \left[\bnabla \cdot \left( |\v|^{q-2} \v \right)\right] d\x \nonumber
    \\
    &\qquad + \nu \int_{\Omega} \v' \cdot \Delta \left( |\v|^{q-2}\v\right) d\x \label{dR}
\end{align}
which, when viewed as a function of the second argument, is a bounded linear functional on both $L^2$ and \H. As an intermediate step, it is convenient to first determine $\nablaL\R_q$, the gradient with respect to the $L^2$ topology. Invoking the Riesz representation theorem \cite{b77}, we have 
\begin{align}
    d\R_q(\v;\v') = \left\langle \v', \nablaL \mathcal{R}_q(\v)\right\rangle_{L^2}
    \label{Riesz}
\end{align}
for all $\v'\in L^2$ which, in view of \eqref{dR}, yields the $L^2$ gradient
\begin{align}
    \nablaL \mathcal{R}_q(\v)&= (q-2) |\v|^{q-4} \v \cdot \big[ \nu \Delta \v + \bnabla \Delta^{-1} \left(\bnabla \v\mathcolon\bnabla \v^T\right)\big] \v \nonumber \\ 
    &\qquad + |\v|^{q-2} \left[\nu \Delta \v + \bnabla \Delta^{-1} \left(\bnabla \v\mathcolon\bnabla \v^T\right)\right]
    \nonumber \\ 
    &\qquad + 2 \partial_i \Delta^{-1} \left[\bnabla \cdot \left( |\v|^{q-2} \v \right)\right] \bnabla \v_i  + \nu \Delta \left( |\v|^{q-2}\v\right).
    \label{gradL2}
\end{align}
Let the Fourier coefficient of the vector field $\v$ be denoted by
\begin{align}
    [ \hatv ]_{\bxi} := \int_{\Omega} \v(\x) e^{- 2\pi i \bxi \cdot \x} \, d\x \in \CC^3,
\end{align}
where $\bxi \in \ZZ^3$ and $i := \sqrt{-1}$. Invoking the Riesz representation theorem again and also using Plancherel's theorem, we obtain
\begin{align}
    \sum_{\bxi \in \ZZ^3}  [\hatvp]_{\bxi}^* \cdot \left[\widehat{\nablaL \mathcal{R}_q}(\v) \right]_{\bxi} &= \left\langle \v', \nablaL \mathcal{R}_q(\v)\right\rangle_{L^2}
    \\
    &= d\R_q(\v;\v')
    \\
    &=\left\langle \v', \nablaH\mathcal{R}_q(\v)\right\rangle_{\H}
    \\
    &= \sum_{l=0,\frac{3}{2}-\frac{1}{q}} \int_{\Omega} \bnabla^l \v'\mathcolon \bnabla^l \nablaH \mathcal{R}_q(\v)\ d\x
    \\
    &= \sum_{l=0,\frac{3}{2}-\frac{1}{q}} \sum_{\bxi \in \ZZ^3} \left((i\boldsymbol{\xi})^l [\hatvp]_{\bxi}\right)^* \mathcolon \, (i\boldsymbol{\xi})^l \left[\widehat{\nablaH \mathcal{R}_q}(\v) \right]_{\bxi}
    \\
    &= \sum_{\bxi \in \ZZ^3}  [\hatvp]_{\bxi}^*  \cdot \left[\left(1+|\boldsymbol{\xi}|^{3-\frac{2}{q}}\right)  \widehat{\nablaH \mathcal{R}_q}(\v)\right]_{\bxi} \label{fracLap}
\end{align}
for all $\v'\in \text{\H{}}$, where ${}^*$ denotes complex conjugation, which implies
\begin{align}
    \left[\widehat{\nablaH \mathcal{R}_q}(\v)\right]_{\bxi} &= \frac{1}{1+|\boldsymbol{\xi}|^{3-\frac{2}{q}}} \left[\widehat{\nablaL \mathcal{R}_q}(\v)\right]_{\bxi}, \qquad \bxi \in \ZZ^3.
    \label{LtoHgradient}
\end{align}
Thus, the Sobolev gradient $\nablaH\R_q \in \text{\H{}}$ can be viewed as the solution to the elliptic boundary-value problem
\begin{equation}
    \left(1 - \Delta^{\frac{3}{2} - \frac{1}{q}}\right) \nablaH\R_q = \nablaL\R_q \qquad \text{in} \ \Omega
\end{equation}
subject to the periodic boundary conditions.

\subsection{Projection Operator}
The operator representing the projection onto the space $\TukXB$ tangent to $X_B$ at $\u_k$ is given by 
\begin{align}
    \mathbb{P}_k: \H \to \mathcal{T}_{\u_k} X_B, \quad \mathbb{P}_{k} \v \defeq \mathbb{P}_{\#}\mathbb{P}_{\sigma}\left( \v - \langle \v, \n_k\rangle_{\H} \n_k \right),
    \label{projectionDef}
\end{align}
where
\begin{align}
    \mathbb{P}_{\sigma}\v \defeq \v - \bnabla \Delta^{-1} \bnabla \cdot \v
\end{align}
is the Leray projection, whereas
\begin{align}
    \mathbb{P}_{\#}\v \defeq \v- \int_{\Omega} \v\ d\x
\end{align}
is the projection onto the subspace of zero-mean functions and $\n_{k}$ is the unit element normal to $X_B$ in \H.

To determine $\n_k$, we define the function $F(\v)\defeq\ \|\v\|_q -B$ such that we have $F(\u) = 0$ for all $\u \in X_B$. Similarly to what we did in Section \ref{subsectionHilberSobolevGradient}, we calculate the gradient of $F$ in \H. Computing the Gateaux differential of $F$, we obtain
\begin{align}
    dF(\v;\v') & = \lim_{\epsilon\to 0}\epsilon^{-1}\left(\|\v+\epsilon \v'\|_q-\|\v\|_q\right)
    = \|\v\|_q^{1-q}\int_\Omega |\v|^{q-2} \v\cdot \v' \ d\x \\
    & = \left\langle \v', \nablaL F(\v) \right\rangle_{L^2} = \left\langle \v', \nablaH F(\v) \right\rangle_{\H},
\end{align}
where we have invoked the Riesz representation theorem, respectively, in $L^2$ and in \H. By following the same steps as in \eqref{Riesz}--\eqref{LtoHgradient}, this allows us to identify the normal elements defined with respect to these two topologies as
\begin{align}
    \nablaL F(\v) & = \|\v\|_q^{1-q}|\v|^{q-2}\v
\end{align}
and
\begin{align}
    \left[\widehat{\nablaH F(\v)}\right]_{\bxi} = \frac{1}{1+|\boldsymbol{\xi}|^{3-\frac{2}{q}}} \left[\widehat{\nablaL F(\v)} \right]_{\bxi}.
\end{align}
Finally, the required unit normal is obtained by performing normalization to give
\begin{align}
    \n_k = \frac{\nablaH F(\u_k)}{\big\|\nablaH F(\u_k)\big\|_{\H}}.
\label{nk}
\end{align}

\subsection{Retraction}
Given that that the norm constraint is expressed in terms of a homogeneous function, the retraction reduces to a simple normalization 
\begin{align}
    \Retraction: \H\to X_B, \qquad 
    \Retraction(\u_k) := B \frac{\u_k\ }{\|\u_k\|_q}.
    \label{retractionDefinition}
\end{align}

\subsection{Optimal Step Size}
The step size $\tau_k$ is found optimally by solving the arc-search problem
\begin{align}
    \tau_k = \argmax_{\tau>0}\R_q\left(\u_{k+1}\right) = \argmax_{\tau>0}\R_q\left(B\frac{\u_k + \tau_k \d_k\ }{\|\u_k + \tau_k \d_k\|_q} \right),
    \label{tauk}
\end{align}
where the objective functional \eqref{defRq} is maximized along an arc on the manifold $X_B$ originating at $\u_k$ in the direction $\d_k$. This is done using a straightforward generalization of Brent's algorithm \cite[Chapter 10]{pressTeukolskyVetterlingFlannery}.

%\subsubsection{optimization procedure}
\subsection{Vector Transport}
For $\u\in X_B$, the vector transport $\Gamma_{\w}(\v)$ describes how the vector $\v \in \TuXB$ is transported in the direction $\w\in \TuXB$. In general, it is given by \cite[Section 8.1]{absilMahonySepulchre}
\begin{align}
    \Gamma: \TXB \oplus \TXB \to \TXB,
    \qquad 
    \Gamma_{\w}(\v) \defeq & \left.\frac{d}{dt} \Retraction\left( \u+\boldsymbol{w}+t\v \right)\right\vert_{t=0},
\end{align}
where $\TXB := \bigcup_{\u \in X_B} \TuXB$ is the tangent bundle on $X_B$. With the retraction given in \eqref{retractionDefinition}, it follows that
\begin{align}
    \Gamma_{\w}(\v) 
    %\defeq & \left.\frac{d}{dt} \RR\left( \u+\boldsymbol{w}+t\v \right)\right\vert_{t=0} \\
    =& \frac{B\v}{\|\u+\boldsymbol{w}\|_q} - \frac{B(\u+\boldsymbol{w})}{ \|\u+\boldsymbol{w}\|_q^{q+1}} \int_\Omega |\u+\boldsymbol{w}|^{q-2}(\u+\boldsymbol{w})\cdot \v \ d\x
    \label{VT}
\end{align}
for all $\v, \w \in \TuXB$.
%Here, we want to transport the previous step direction $\d_{k-1}$ along the previous step $\tau_{k-1} \d_{k-1}$, which yields
%\begin{align}
%    \Gamma_{\tau_{k-1}\d_{k-1}}&(\d_{k-1})
%    \\
%    &= \frac{B\v}{\|\u_{k-1}+\tau_{k-1}\d_{k-1}\|_q} - \frac{B(\u_{k-1}+\tau_{k-1}\d_{k-1})}{ \|\u_{k-1}+\tau_{k-1}\d_{k-1}\|_q^{q+1}}
%    \\
%    &\qquad\qquad\cdot\int_\Omega |\u_{k-1}+\tau_{k-1}\d_{k-1}|^{q-2}(\u_{k-1}+\tau_{k-1}\d_{k-1})\cdot \v \ d\x
%\end{align}
\subsection{Momentum Coefficient}
We chose the Polak-Ribiere version of the nonlinear conjugate-gradient method with the momentum coefficient given by \cite{nocedalWright}
\begin{align}
    %\beta_k &= \frac{\left\langle\mathbb{P}_{k} \left(\nablaH\R_q(\u_k)\right),\mathbb{P}_{k} \left(\nablaH\R_q(\u_k)\right)-\Gamma_{\tau_{k-1}\d_{k-1}}\left(\mathbb{P}_{k-1} \left(\nablaH\R_q(\u_{k-1})\right)\right)\right\rangle_\H}{\left\|\mathbb{P}_{k-1} \left(\nablaH\R_q(\u_{k-1})\right)\right\|_{\H}^2}
    \beta_k &= \frac{\left\langle(\mathbb{P}\nabla\R)_{k},(\mathbb{P}\nabla\R)_{k} -\Gamma_{\tau_{k-1}\d_{k-1}}\left((\mathbb{P}\nabla\R)_{{k-1}}\right)\right\rangle_{\H}}{\left\|(\mathbb{P}\nabla\R)_{k-1}\right\|_{\H}^2},
    \label{beta_def}
\end{align}
where $(\mathbb{P}\nabla\R)_{k} := \mathbb{P}_{k} \left(\nablaH\R_q(\u_k)\right)$.

\subsection{Iteration Scheme}
\begin{figure}[ht]
    \centering
    \resizebox{\linewidth}{!}{\definecolor{mcmaster_maroon}{RGB}{122, 0, 60}%
\definecolor{mcmaster_gold}{RGB}{253, 191, 87}%
\definecolor{mcmaster_grey}{RGB}{94, 106, 113}%

\colorlet{colorA}{mcmaster_maroon}%
\colorlet{colorB}{mcmaster_gold}%
\colorlet{colorC}{mcmaster_grey}%

\colorlet{colorA}{paraviewBlue}%
\colorlet{colorB}{paraviewRed}%
\colorlet{colorC}{matplotlibGreen}%

\def \r {15}%
\def \circlePhiStart {62.5}%
\def \circlePhiEnd {180-\circlePhiStart}%

\def \rUkbasePoint {7*\r/8}%
\def \rUkminusOnebasePoint {29*\r/32}%

\def \phiUminusOne {180-\circlePhiStart-7.5}%
\def \xUminusOne {(\r*cos(\phiUminusOne))}%
\def \yUminusOne {(\r*sin(\phiUminusOne))}%

\def \xDeltaGrad {2}%
\def \yDeltaGrad {4}%

\def \innerMinusOne {(sin(\phiUminusOne)*(\xDeltaGrad) -cos(\phiUminusOne)*(\yDeltaGrad))}%

\def \xProjDelta {(\innerMinusOne*sin(\phiUminusOne))}%
\def \yProjDelta {(-\innerMinusOne*cos(\phiUminusOne))}%

\def \tauLength {2.97}%

\def \xTauProjDelta {(\tauLength*\innerMinusOne*sin(\phiUminusOne))}%
\def \yTauProjDelta {(-\tauLength*\innerMinusOne*cos(\phiUminusOne))}%

\def \phiU {(90-atan((\xUminusOne+\xTauProjDelta)/(\yUminusOne+\yTauProjDelta)))}%

\def \xU {(\r*cos(\phiU))}%
\def \yU {(\r*sin(\phiU))}%

\def \inner {(sin(\phiU)*(\xProjDelta) -cos(\phiU)*(\yProjDelta))}%

\def \xGammaDelta {(\inner*sin(\phiU))}%
\def \yGammaDelta {(-\inner*cos(\phiU))}%

\def \dashedOverShoot {0.25}%

\def \xDeltaGradUK {-0.375}%
\def \yDeltaGradUK {1.5}%

\def \innerGradUK {(sin(\phiU)*(\xDeltaGradUK) -cos(\phiU)*(\yDeltaGradUK))}%

\def \xDeltaProjGradUK {(\innerGradUK*sin(\phiU))}%
\def \yDeltaProjGradUK {(-\innerGradUK*cos(\phiU))}%

\def \phiProjGradUK {(90-atan((\xU+\xDeltaProjGradUK)/(\yU+\yDeltaProjGradUK)))}%

\def \deltaHeight {\height*1/24}%

\begin{tikzpicture}[line cap=round, line width = 0.8pt]

    % padding 
    %\node at (-5,7) {};
    %\node at (5,0) {};
    
    \draw[samples = 1000, domain=\circlePhiEnd:\circlePhiStart] plot ({\r*cos(\x)}, {\r*sin(\x)}) node[below] {$X_{B}$};

    \draw[->, colorC]  ({\rUkminusOnebasePoint*cos(\phiUminusOne)},{\rUkminusOnebasePoint*sin(\phiUminusOne)}) -- ({\xUminusOne},{\yUminusOne}) node[below right] {$\u_{k-1}$};

    \draw[->, colorA, line width = 1.4pt]  ({\xUminusOne},{\yUminusOne}) -- +({\xProjDelta},{\yProjDelta}) node[above left] {$\d_{k-1}$};
    
    \draw[->, colorC]  ({\xUminusOne},{\yUminusOne}) -- +({\xTauProjDelta},{\yTauProjDelta}) node[above left] {$\tau_{k-1}\d_{k-1}$};

    \draw[dashed, line width = 0.4] ({\xU},{\yU}) -- ({\xUminusOne+\xTauProjDelta},{\yUminusOne+\yTauProjDelta});

    \draw[->, colorC]  ({\rUkbasePoint*cos(\phiU)},{\rUkbasePoint*sin(\phiU)}) -- ({\xU},{\yU}) node[below left] {$\u_{k}$};

    \draw[dashed, line width = 0.4] ({\xU+\xGammaDelta},{\yU+\yGammaDelta}) -- ({\xU+\xProjDelta},{\yU+\yProjDelta});
    
    \draw[dashed, line width = 0.4] ({\xU+\xDeltaProjGradUK-\dashedOverShoot*sin(\phiU)},{\yU+\yDeltaProjGradUK+\dashedOverShoot*cos(\phiU)}) -- ({\xU+\xGammaDelta+\dashedOverShoot*sin(\phiU)},{\yU+\yGammaDelta-\dashedOverShoot*cos(\phiU)});

    \draw[dashed, line width = 0.4] ({\xU+\xDeltaProjGradUK},{\yU+\yDeltaProjGradUK}) -- ({\xU+\xDeltaGradUK},{\yU+\yDeltaGradUK});

    \draw[->, colorA]  ({\xU},{\yU}) -- +({\xProjDelta},{\yProjDelta}) node[right, above] {$\d_{k-1}$};

    \draw[->, colorA]  ({\xU},{\yU}) -- +({\xGammaDelta},{\yGammaDelta}) node[right, yshift = 5] {$\Gamma (\d_{k-1})$};

    \draw[->, colorB]  ({\xU},{\yU}) -- +({\xDeltaGradUK},{\yDeltaGradUK}) node[above] {$\nablaH\R(\u_{k})$};

    \draw[->, colorB]  ({\xU},{\yU}) -- +({\xDeltaProjGradUK},{\yDeltaProjGradUK}) node[above left, xshift = 10] {$\mathbb{P}_k \nablaH\R(\u_{k})$};
    
    %\draw[->, thin, yellow]  (0,0) -- ({2*\r*cos(\phiProjGradUK)},{2*\r*sin(\phiProjGradUK)}) node[below left] {$test nonsense$};

\end{tikzpicture}}
    \caption{Schematic illustration of a single iteration of the Riemannian conjugate gradient method described in \eqref{RCG1}--\eqref{RCG2}, based on \cite[Fig. 1]{zhaoProtas2023}.}
    \label{fig:RCG}
\end{figure}
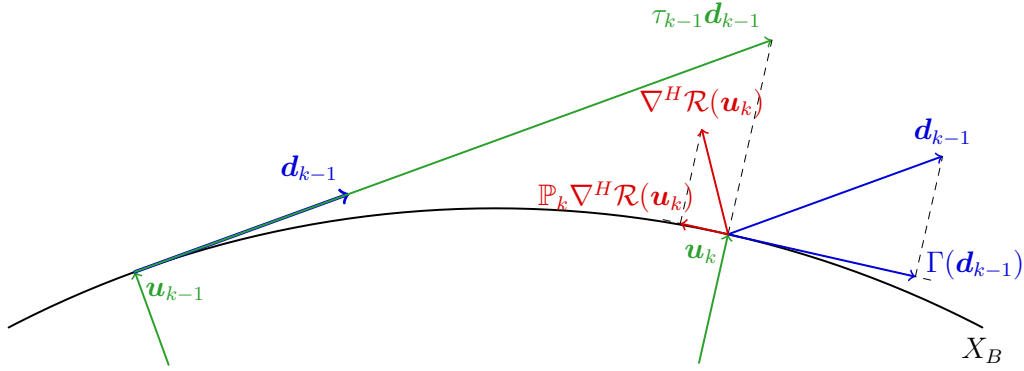

The full iteration scheme as derived in the previous subsections is given by
\begin{align}
    \u_{k+1} &=  B \frac{\u_k + \tau_k \d_k}{\|\u_k + \tau_k \d_k\|_q} \qquad k =0,1,2, \dots \label{RCG1}
    \\
    \d_k &= \mathbb{P}_{k} \left(\nablaH\R_q(\u_k)\right) + \beta_k\Gamma_{\tau_{k-1}\d_{k-1}}(\d_{k-1})\label{RCG2}
\end{align}
with $\d_0=\0$, where
\begin{itemize}
    \item $\nablaH\R_q(\u_k)$ is defined via \eqref{LtoHgradient} and \eqref{gradL2},
    \item $\mathbb{P}_{k}$ is defined in \eqref{projectionDef},
    \item $\tau_k$ is found by solving \eqref{tauk},
    \item $\beta_k$ is defined in \eqref{beta_def},
    \item $\Gamma_{\tau_{k-1}\d_{k-1}}(\d_{k-1})$ is defined in \eqref{VT}.
\end{itemize}
The algorithm represented by \eqref{RCG1}--\eqref{RCG2} is schematically illustrated in Figure \ref{fig:RCG}.

%\subsubsection{computation of branches via continuation and initialization}
\subsection{Computation of Branches via Continuation}
In order to test the sharpness of estimate \eqref{estimate}, we aim to construct solutions of Problem \ref{hilbertProblem} for different $q$ and increasing values of the constraint parameter $B$. For each value of $B$ the approach described above requires an initial guess $\u_0$. For the smallest considered value of $B$ the initial guess $\u_0=\uInit$ is chosen as the ABC flow \eqref{ABC} which is the local maximizer of the problem in the small-data limit corresponding to the maximum value of the objective functional \eqref{defRq}, cf.~Section \ref{sectionSmallDataLimit}. Then, for increasing values of this parameter $B$, we use a continuation approach where the optimal solution $\tuB$ serves as the initial guess $\u_0^{B(1+\delta B)}$ to compute $\widetilde{\u}^{B(1 + \delta B)}$, i.e., the solution corresponding to the next, incremented, value of the constraint parameter. This approach is summarized as Algorithm \ref{alg1}.

\begin{algorithm}[htbp]
\DontPrintSemicolon
\hrule
\BlankLine
\KwIn{\\
\Indp
$q$ --- Exponent of the Lebesgue space of interest\\
$\uInit$ --- Initial guess corresponding to the small-data limit (see Section \ref{sectionSmallDataLimit})\\
$B^{\text{init}}$ --- Initial value for the constraint parameter\\
$B^{\text{max}}$ --- Maximum value for the constraint parameter\\
$\delta B$ --- Increment of the constraint parameter \\
$\epsilon$ --- Convergence tolerance\\
\Indm
}
\BlankLine
\hrule
\BlankLine
\KwOut{\\
\Indp
Extreme states $\tuB$ and corresponding values of the objective functional $\R(\tuB)$ for $B^{\text{init}}\leq B \leq B^{\text{max}}$\\
\Indm
}
\BlankLine
\hrule
\BlankLine

\SetKwFunction{main}{MainBranch}
\SetKwFunction{solveHilbertProb}{SolveProblem\ref{hilbertProblem}}
\SetKwProg{Fn}{Function}{:}{}
\Fn{\main{}}{
    $B\leftarrow B^{\text{init}}$\\
    $\u_0\leftarrow\uInit$\\
    \Repeat{$B\geq B^{\text{max}}$}{
        $\tuB \leftarrow\ $\solveHilbertProb{$q$,$B$,$\u_0$,$\epsilon$}\\
        Evaluate $\R(\tuB)$ according to \eqref{defRq}\\
        $\u_0\leftarrow\tuB$\tcp*{Branch continuation with prev max}
        $B\leftarrow B(1+\delta B)$\tcp*{Increase constraint value}
    }
}
\BlankLine
\SetKwProg{Fn}{Function}{:}{}
\Fn{\solveHilbertProb{$q$,$B$,$\u_0$,$\epsilon$}}{
    $k=0$\\
    $\Gamma_{\tau_{-1}\d_{-1}}(\d_{-1})=\0$\\
    \Repeat{$\left|\R(\u_{k+1})-\R(\u_k)\right| / |\R(\u_k)| <\epsilon$}{%
        Compute the $L^2$ gradient $\nablaL \R(\u_k)$ according to \eqref{gradL2}\\
        Compute the \H{} gradient $\nablaH \R(\u_k)$ according to \eqref{LtoHgradient}\\
        Compute the step direction $\d_k$ according to \eqref{RCG2}\\
        Find step size $\tau_k$ by solving \eqref{tauk}\\
        Compute $\u_{k+1}$ by \eqref{RCG1}\\
        Compute the vector transport $\Gamma_{\tau_{k}\d_{k}}(\d_{k})$ by $\eqref{VT}$ for the next iteration\\
        $k\leftarrow k+1$\\
    }
    \KwRet $\u_{k+1}$
}
\BlankLine
\caption{Continuation approach used to compute branches of local maximziers $\tuB$ for increasing values of the constraint parameter $B$.}
\label{alg1}% label after caption!!!
\end{algorithm}

%\subsubsection{numerics (pseudo-spectral discretization, dealiasign, regularity)}

\subsection{Spatial Discretization}
Here we provide information about the spatial discretization used in different elements of Algorithm \ref{alg1} (the objective functional \eqref{defRq}, gradient \eqref{gradL2}, etc.). We use a standard pseudospectral Fourier-Galerkin method with products evaluated in the physical space \cite{canutoHussainiQuarteroniZang}. As regards dealiasing, several expressions that need to be evaluated involve products of degree higher than two and in principle it is possible to adapt the 3/2 rule to these cases. However, from the point of view of the computational efficiency, it is better to factor all products into quadratic products which are then evaluated using the standard 3/2 rule. The algorithm is implemented in {\tt Fortran 90}  with {\tt MPI}  and Fourier transforms are performed using the library {\tt FFTW3}. 

To improve the convergence of the iterations \eqref{RCG1}--\eqref{RCG2}, the inner product \eqref{Hs} is redefined as 
\begin{align}
    \langle \boldsymbol{f},\boldsymbol{g}\rangle_{H^s} &\defeq  \int_{\Omega}  \boldsymbol{f}(\x)\cdot \boldsymbol{g}(\x) + \ell^{2s}\bnabla^{s} \boldsymbol{f}(\x)\mathcolon\bnabla^{s} \boldsymbol{g}(\x)\, d\x,
    \label{Hsl}
\end{align}
where $\ell \in (0,\infty)$ is an adjustable parameter. Clearly, \eqref{Hsl} is an equivalent inner product to \eqref{Hs}. As shown in \cite{protasBewleyHagen2004}, computation of Hilbert-Sobolev gradients based on $L^2$ gradients can be viewed as low-pass filtering of the latter, cf.~\eqref{LtoHgradient}, where $\ell$ acts as the cut-off length scale. Adjusting the value of this parameter helps improve the convergence of iterations in Algorithm \ref{alg1} and, based on a number of tests, the value we use is $\ell = 0.1$. In addition, the momentum term $\beta_k$, cf.~\eqref{RCG2}, is reset to zero every 25 iterations which also improves convergence.

In the course of iterations, we carefully monitor the Fourier spectrum of the solutions $\u_k$ to make sure they are always well resolved and the resolution is refined (by doubling the number $N$ of grid points in each direction) when the resolution limit is approached. The largest resolution used in the results reported in the next section is $N^3 = 1024^3$.

\subsection{Physical and Numerical Parameters}
The physical and numerical parameters that were mostly used in this investigation are given in Table \ref{tab:numericalParameters}.
\begin{table}[ht]
    \centering
    \begin{tabular}{c|c|c|c}
        Parameter & Variable & Value & Remark \\
        \hline
        Viscosity & $\nu$ & $1$ & see \eqref{nse} \\
        Resolution & $N$ & $512^3$, $1024^3$ & see Fig.~\ref{fig:branches} for refinement\\
        Convergence tolerance & $\epsilon$ & $10^{-5}$, $10^{-4}$ & see Algorithm \ref{alg1} \\
        Initial constraint size & $B^{\text{init}}$ & $1$ & see Algorithm \ref{alg1} \\
        Constraint increment & $\delta B$ & $10^{1/8}-1$ & see Algorithm \ref{alg1} \\
        %Maximal constraint size & $\dots$  & $B^{\text{max}}$ in Algorithm \ref{alg1} \\
        Sobolev parameter & $\ell$ & $0.1$ & see \eqref{Hsl}
    \end{tabular}
    \caption{Values of the parameters used in the computations reported in Section \ref{sectionResults}.}
    \label{tab:numericalParameters}
\end{table}

%\subsubsection{gradient - Hilbert vs Banach space formulation}
\subsection{Banach-Sobolev Gradients}
Finally, we offer some comments about computation of gradients that would be required to solve Problem \ref{banachProblem}. The main difficulty is that this problem is formulated in the space $W^{1,\frac{3q}{q+1}}$, which is not endowed with an inner product for $q>3$ and therefore one cannot make use of the Riesz theorem, as was done above where we constructed the Hilbert-Sobolev gradients, cf.~\eqref{Riesz}. To get around this difficulty, the Banach gradient can be defined as a solution of the following variational optimization problem \cite{neuberger,protas2008,ramirez2025}
\begin{align}
    \nablaB \mathcal{R}_q(\u) = \argmax_{\substack{\v\in W\\\|\v\|_W=\|d\mathcal{R}_q(\u;\cdot)\|_{W^\star}}} d\mathcal{R}_q(\u;\v) = \argmax_{\substack{\v\in W\\ \|\v\|_W=\|\nablaL\mathcal{R}_q\|_{W^\star}}} \langle \nablaL \mathcal{R}_q, \v\rangle_{L^2},\quad
    \label{dRB}
\end{align}
where $W := W^{1,\frac{3q}{q+1}}$ and $W^\star$ is its dual. It is motivated by the observation that the gradient is the element maximizing the differential $d\R_q(\u;\cdot)$ and reduces to the classical formulation described above when $W$ is a Hilbert space \cite{neuberger}. Problem \eqref{dRB} can be converted to an unconstrained formulation by introducing the Lagrange multipliers $\rho \in L^\frac{3q}{2q-1}$ and $\lambda \in \RR$
\begin{align}
    \nablaB \mathcal{R}_q(\u) = \min_{\rho \in L^\frac{3q}{2q-1}, \lambda \in \mathbb{R}_q}\argmax_{\v\in W} &\int_\Omega \bigg[\nablaL\mathcal{R}\cdot \v + \bnabla \cdot \v \rho +\notag
    \\
    & \frac{\lambda(q+1)}{3q} \left( |\v|^{\frac{3q}{q+1}}+ |\bnabla \v|^{\frac{3q}{q+1}} -  \frac{1}{|\Omega|}\right)\bigg] d\x.
\end{align}
and the corresponding Euler-Lagrange equations are 
\begin{gather}
    |\nablaB\mathcal{R}_q|^{\frac{q-2}{q+1}}\nablaB\mathcal{R}_q - \sum_{i=i}^3 \partial_{x_i} \left(\left|\bnabla \left(\nablaB\mathcal{R}_q\right)\right|^{\frac{q-2}{q+1}} \partial_{x_i} \nablaB\mathcal{R}_q\right) = \frac{1}{\lambda} \bnabla\rho - \frac{1}{\lambda} \nablaL\mathcal{R}_q, \label{EL1_Banach} \\
    \Delta \rho = \lambda \bnabla\cdot\left[|\nabla^W\mathcal{R}_q|^{\frac{q-2}{q+1}}\nabla^W\mathcal{R}_q - \partial_{x_i} \left(\left|\bnabla \left(\nabla^W\mathcal{R}_q\right)\right|^{\frac{q-2}{q+1}} \partial_{x_i} \nabla^W\mathcal{R}_q\right)\right] \label{EL2_Banach}
\end{gather}
with $\lambda$ determined by the normalization condition $\|\nablaB\R_q\|_W = 1$. The main difference with respect to how the Hilbert-Sobolev gradients are computed (see Subsection \ref{subsectionHilberSobolevGradient}) is that now the relation between the $L^2$ and the Hilbert-Banach gradients $\nablaL\R_q$ and $\nablaB\R_q$ is no longer linear. More specifically, equation \eqref{EL1_Banach} is a nonlinear boundary-value problem involving the $p$-Laplacian with $p = (q - 2) / (q + 1)$. Since the numerical solution of such a problem typically requires a variant of Newton's method where the Jacobian needs to be evaluated, this problem is intractable in the present setting given the required numerical resolutions.

\section{Results}
\label{sectionResults}
%\subsubsection{diagnostics for a representative optimization - \texorpdfstring{$R(u^n) \text{ vs }n$}{R of un vs n}, spectra}
In this section we present the results obtained by solving Problem \ref{hilbertProblem} for $q = 4,5,6,9$ and increasing values of the constraint parameter $B$. Families of local maximizers $\tuB$ parameterized by $B$ and computed with Algorithm \ref{alg1} are referred to as "branches". In the computations reported here we set $\nu = 1$ in \eqref{nse}, and use the numerical resolutions of $N^3 = 512^3$ and $N^3 = 1024^3$, respectively, for "small" and "large" values of $B$. Below we first discuss the solution of Problem \ref{hilbertProblem} for representative values of $q$ and $B$ before analyzing the branches obtained for different $q$ in order to draw conclusions about the sharpness of estimate \eqref{estimate}. We also analyze the structure of the local maximizers $\tuB$ and provide some comments about solutions of Problem \ref{hilbertProblem} for $q \rightarrow 3$.

\subsection{Solution of Problem \ref{hilbertProblem} for Representative Values of \texorpdfstring{$q$}{q} and \texorpdfstring{$B$}{B}}
\label{subsectionRepresentativeValues}

\begin{figure}[htp] % allows for here, top or as float page, fixes long caption problem -> uses float page
    %\centering    \includegraphics[width=\linewidth]{pictures/iter+spectra_q5_B019.png}
    
    \begin{subfigure}{\linewidth}
        \centering
        \includegraphics[width=\textwidth]{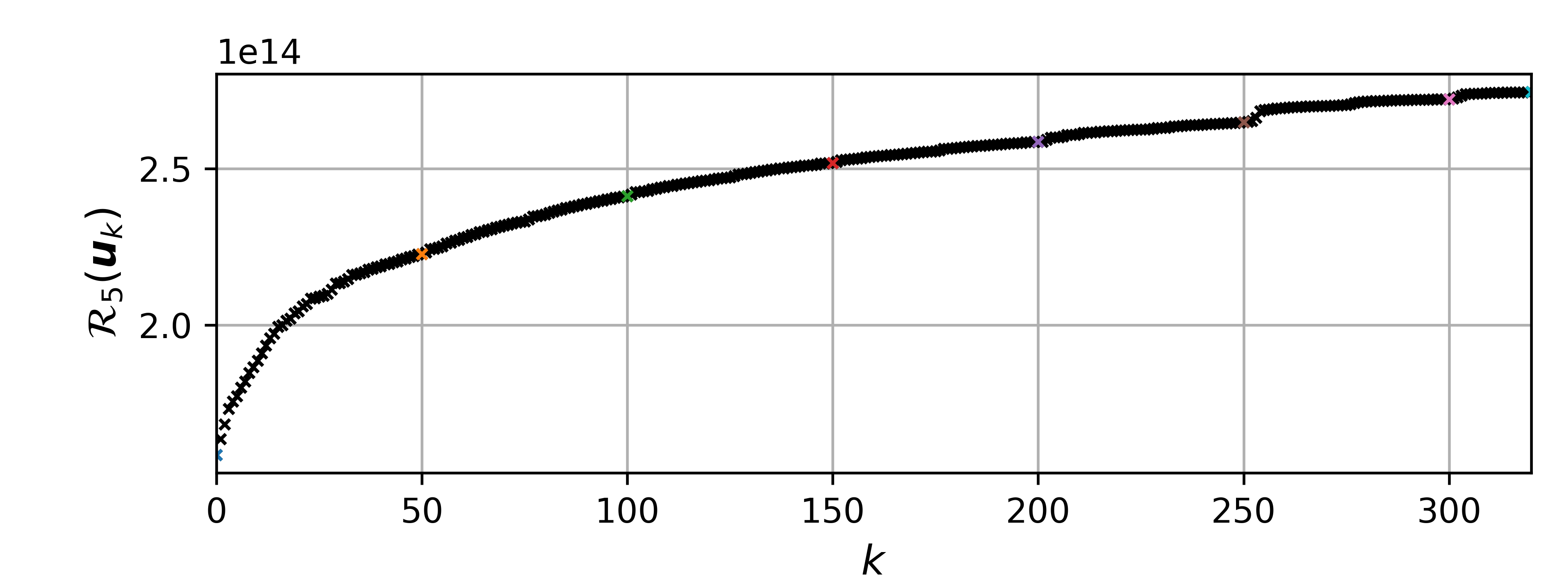}
        \caption{}
        \label{fig:iteration+spectra_aIteration}
    \end{subfigure}
    
    \begin{subfigure}{\linewidth}
        \centering
        \includegraphics[width=\textwidth]{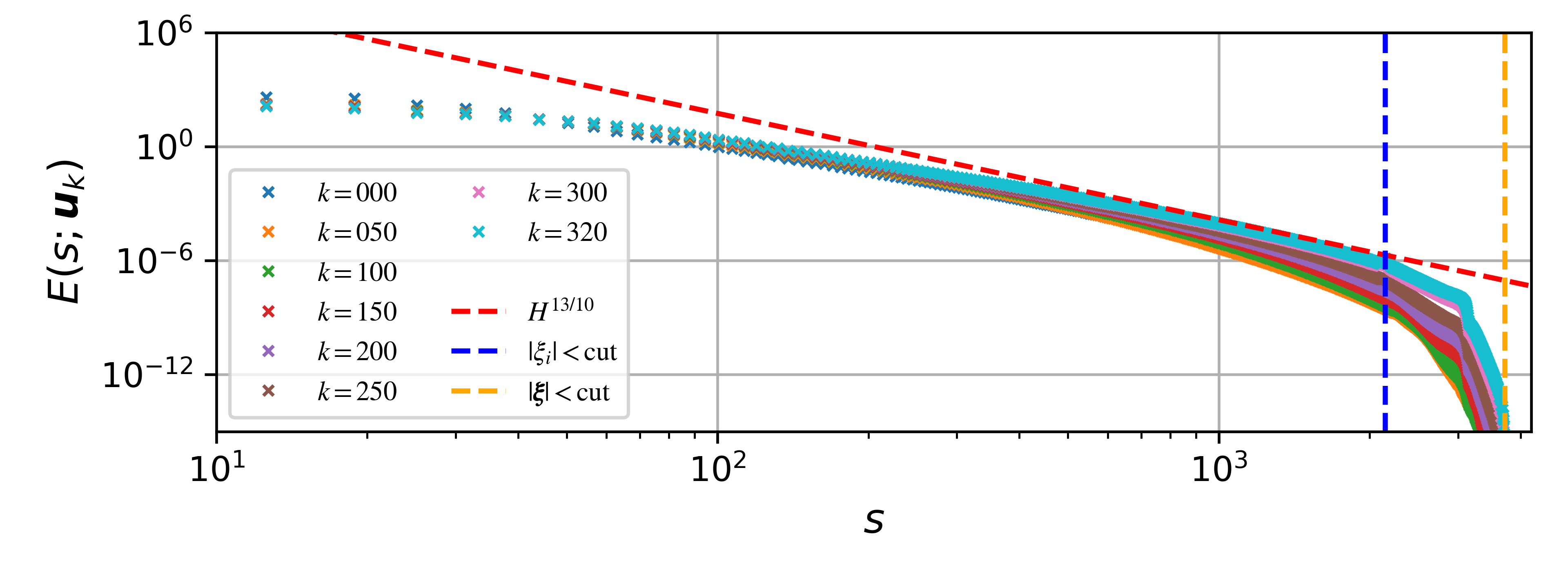}
        \caption{}
        \label{fig:iteration+spectra_bSpectra}
    \end{subfigure}
    \caption{(a) The objective functional $\R_5(\u_k)$ versus iterations $k$ in the solution of Problem \ref{hilbertProblem} with $q=5$ and $B = 10^{9/4}\approx 177.8$. (b) The corresponding energy spectra $E(s,\u_k)$ (see \eqref{spectraDef}) of the approximations $\u_k$ of the maximizer $\tuB$ obtained at different iterations $k$; the slanting {\color{red}red} line represents the relation %$C |\bxi|^{-2\left(\frac{3}{2}-\frac{1}{q}\right) -n} = C s^{-6+\frac{2}{q}} = C s^{-5.6}$
        $C s^{-28/5}$ for some $C > 0$ which is the asymptotic form (as $|\bxi| \rightarrow \infty$) of the energy spectrum for functions immediately outside the edge of regularity in the space $H^{13/10}$; the left {\color{blue}blue} vertical line at $s=\frac{2}{3}\pi N\approx 2145$ represents the wavenumber below which the Fourier components are unaffected by dealiasing whereas the right  {\color{orange}orange} vertical line at $s=\frac{2}{\sqrt{3}}\pi N\approx 3715$ marks the wavenumber above which {\it all} Fourier components are affected (in other words, for the wavenumbers between the two lines only some Fourier components are affected by the filter).}
    \label{fig:iteration+spectra}
\end{figure}

To fix attention, we focus here on the solution of Problem \ref{hilbertProblem} with $q=5$ and $B = 10^{9/4}\approx 177.8$. It is formulated in the space $H^{13/10}$, where $13/10 = 3/2 - 1/q$ %$\frac{13}{10} = \frac{3}{2} - \frac{1}{q}$
with $q = 5$ and numerically solved with a resolution of $N^3=1024^3$. The increase of the objective functional $\R_5(\u_k)$ with iterations $k$ is shown in Figure \ref{fig:iteration+spectra_aIteration} where a monotonic growth is evident. The observed net increase (less than by a factor of two) may not appear large, however, it should be observed that Algorithm \ref{alg1} employs the continuation approach where the previous maximizer $\widetilde{\u}^{B/(1 + \delta B)}$ is used as the initial guess $\u_0$ when solving Problem \ref{hilbertProblem} for the current value of $B$. Therefore, when the increments $\delta B$ are small, the initial value $\R_5(\u_0)$ of the objective functional is already quite close to the maximum $\R_5(\tuB)$. 

%    #yReg = x[:]**(-2.0*s-3.0)*factor       # infty > ||u||_H^s^2 = int (1+|xi|^s)^2 |u(xi)|^2 dxi ~ c int_0^infty |xi|^{2s} |u(|xi|)|^2 |xi|^{n-1} d|xi|
%                                        # borderline for integrant ~ |xi|^{-1} -> log(infty)
%                                        # |u(|xi|)|^2 |xi|^{2s+n-1} ~ -1
%                                        # |u(|xi|)|^2 ~ |xi|^{-2s-n}

In Figure \ref{fig:iteration+spectra_bSpectra} we show the energy spectra
\begin{align}
    E(s; \u_k) :=\sum_{s-2\pi<|\bxi|<s} \left| \left[\widehat \u_k \right]_{\bxi} \right|^2
    \label{spectraDef}
\end{align}%
of the approximations obtained as iterations \eqref{RCG1}--\eqref{RCG2} progress and we conclude that they are all properly resolved. We also observe that while the approximations are smoother at early iterations (which is reflected in the faster decay of $E(s;\u_k)$ as $s \rightarrow \infty$), they approach the "edge of regularity" in the space $H^{13/10}$, as the iterations converge and $\u_k \rightarrow \tuB$. The structure of the optimal fields $\tuB$ in the physical space is discussed in Subsection \ref{subsectionStructure}.

\begin{figure}[ht]
    \centering
    \includegraphics[width=\linewidth]{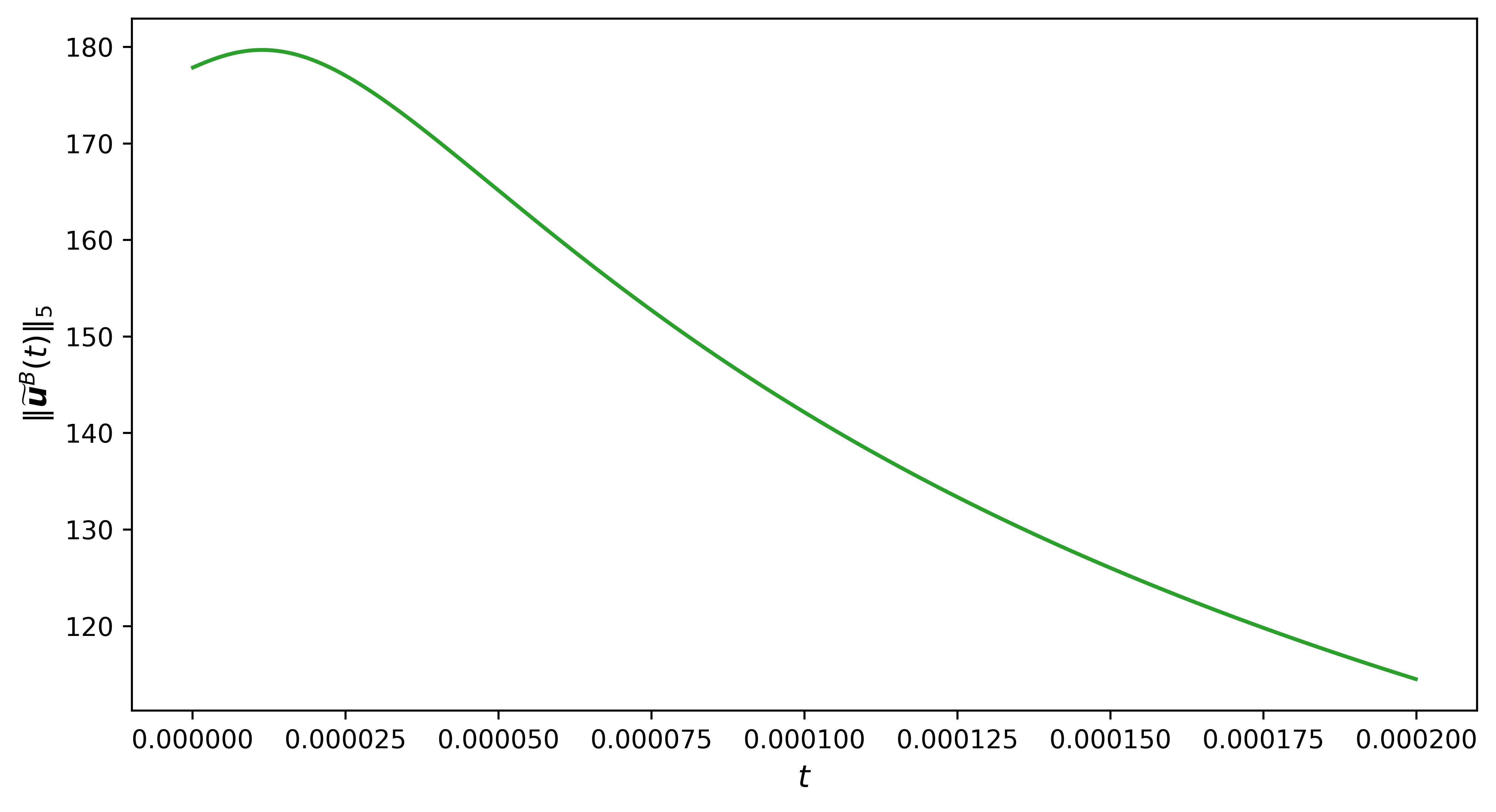}
    \caption{Time evolution of the $L^5$ norm of the solution of system \eqref{nse}--\eqref{incompressible} with the initial condition given by the optimal field $\tuB$ obtained by solving Problem \ref{hilbertProblem} for $q=5$ and $B = 10^{9/4}\approx 177.8$.}
    \label{fig:timeEvolution_q5}
\end{figure}

The time evolution of $\| \u(t) \|_5$ in the solution of the the Navier-Stokes system \eqref{nse}--\eqref{incompressible} with the optimal initial data $\widetilde{\u}^{B}$ discussed above is shown in Figure \ref{fig:timeEvolution_q5}. As expected, since in Problem \ref{hilbertProblem} we seek to maximize the instantaneous growth of $\| \u(t) \|_5$,  this quantity increases only initially for $t = 0^+$ after which its growth is rapidly depleted. This observation is consistent with what was reported in \cite{ayalaProtas2017} for solutions of Problem \ref{pb:maxdEdt}.

%\subsubsection{branches and fits, exponents}

\subsection{Branches of Local Maximizers Obtained for Different \texorpdfstring{$q$}{q}}

\begin{figure}[ht]
    \centering
    \includegraphics[width=\linewidth]{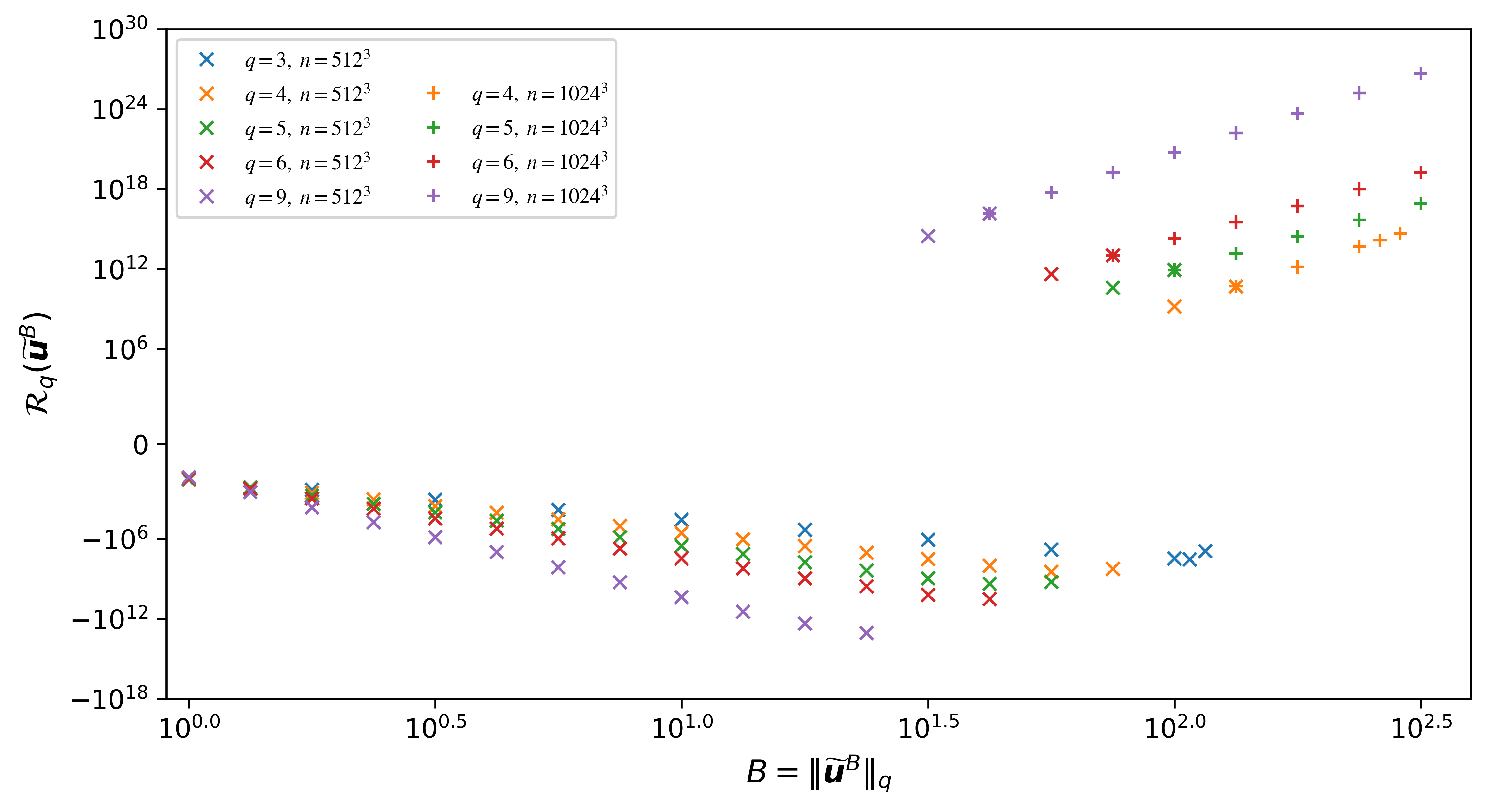}
    \caption{Dependence of the maximum values of the objective functional $\R_q(\tuB)$ on $B$ for $q=$ {\color{matplotlibBlue}(blue)} 3, {\color{matplotlibOrange}(orange)} 4, {\color{matplotlibGreen}(green)} 5, {\color{matplotlibRed}(red)} 6, and {\color{matplotlibPurple}(purple)} 9. The vertical axis is scaled linearly within the interval $[-10^4 : 10^4 ]$ and logarithmically outside.}
    \label{fig:branches}
\end{figure}

\begin{figure}[ht]
    \centering
    \includegraphics[width=\linewidth]{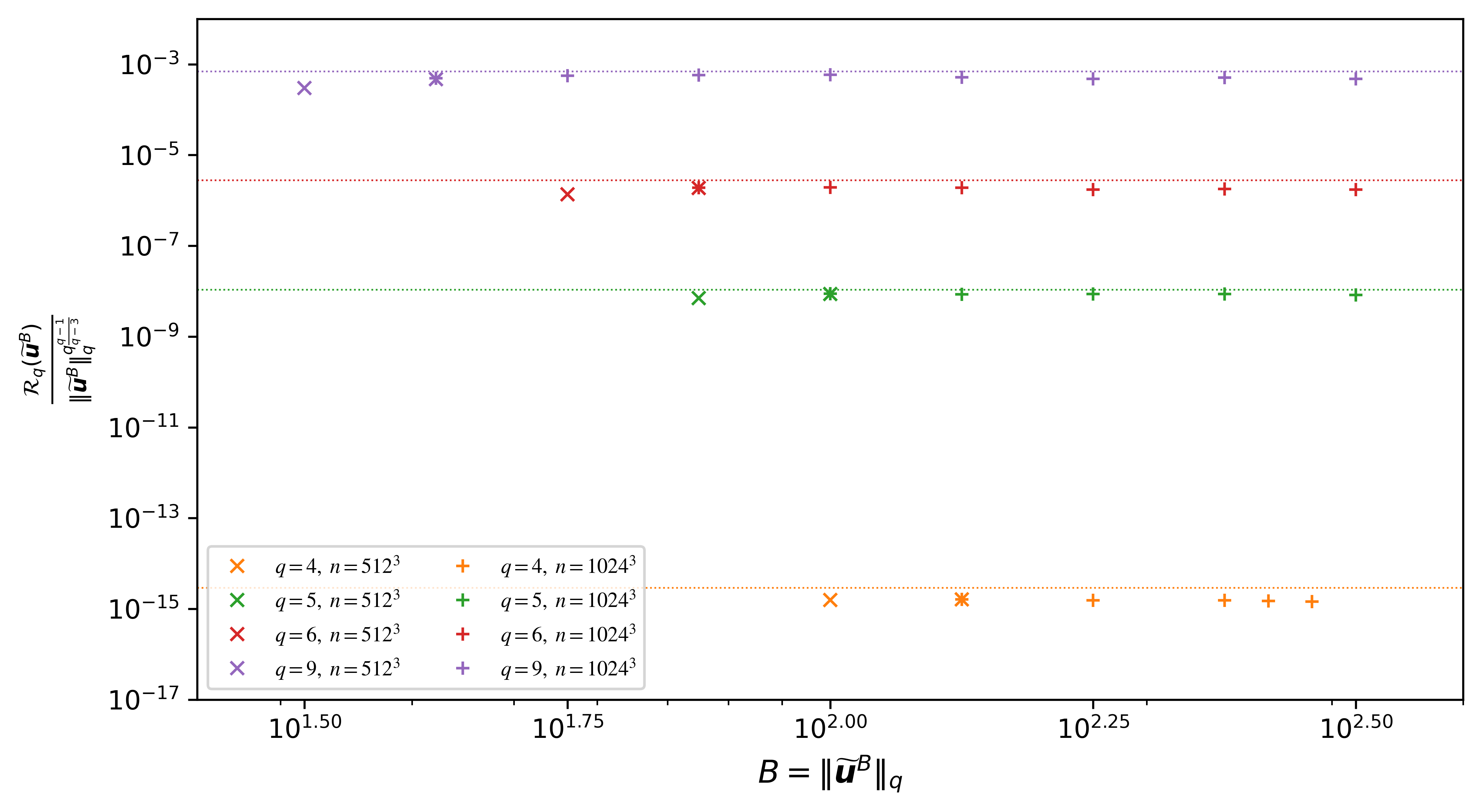}
    \caption{Dependence of the maximum values of the objective functional compensated with the power-law expression on the right-hand side of \eqref{estimate} $\R_q(\tuB) / \| \tuB \|_q^{q(q-1)/(q-3)}$ for $B > \Bchange$ and $q=$ {\color{matplotlibOrange}(orange)} 4, {\color{matplotlibGreen}(green)} 5, {\color{matplotlibRed}(red)} 6, and {\color{matplotlibPurple}(purple)} 9. The elevation of the dashed horizontal lines represents the prefactors $\widetilde C$ from the fitting ansatz \eqref{fit}, cf.~Table \ref{tab:exponents}. All markers lie below their respective constant lines since the fitted exponents are slightly smaller than $q\frac{q-1}{q-3}$.}
    \label{fig:compensated_scaling}
\end{figure}

We now move on to discuss the global picture presented by the branches obtained with Algorithm \ref{alg1} for $q = 3, 4, 5, 6, 9$. Those results are shown in Figure \ref{fig:branches} where two distinct regimes are evident. For each examined value of $q$, except for $q = 3$, there exists a $\Bchange = \Bchange(q) > 0$, such that $\R_q(\widetilde{\u}^{\Bchange}) = 0$ and $\R_q(\tuB) < 0$ for $B \in (0,\Bchange)$, meaning that the rate of growth of the $L^q$ norm is in fact negative in the small-data regime. This is because in this regime the dissipative viscous term in the objective functional \eqref{defRq}, which is of order $q$ in $\u$, dominates the pressure term, which is of order $q+1$ in $\u$. This observation is also consistent with our analysis of solutions to Problem \ref{hilbertProblem} in the limit $B \rightarrow 0$ in Section \ref{sectionSmallDataLimit}, cf.~\eqref{objectiveExpansion2}. On the other hand, for $B > \Bchange$ the pressure term dominates and $\R_q(\tuB) > 0$. Given the range of $B$ shown in Figure \ref{fig:branches}, which spans nearly three orders of magnitude, the transition between these two regimes occurs relatively rapidly. For each value of $q$ the branches are terminated when the maximizers $\tuB$ can no longer be accurately approximated with the resolution $N^3 = 1024^3$. In Figure \ref{fig:branches} we see that as $B \rightarrow \infty$ the value of the objective functional exhibits a power-law behavior in $B$ which can be represented by the ansatz
\begin{equation}
\R_q(\tuB) =\tC \| \tuB \|_q^{\talpha}, \qquad \tC, \ \talpha > 0.
\label{fit}
\end{equation}
In order to quantify this observation and make a connection with estimate \eqref{estimate}, we plot the maximum values of the objective functional  compensated with the power-law in the upper bound in \eqref{estimate}, $\R_q(\tuB) / \| \tuB \|_q^{q(q-1)/(q-3)}$, as a function of $B$ in Figure \ref{fig:compensated_scaling}. The fact that this data falls on horizontal lines suggests that the upper bound in \eqref{estimate} is saturated in terms of the exponent by $\R_q(\tuB)$ as $B \rightarrow \infty$ with the prefactor $\tC$ represented by the vertical elevation of the lines. The exponents on the right-hand side in bound \eqref{estimate} and the exponents $\talpha$ in the ansatz above found by performing least-squares fits to the data in Figure \ref{fig:branches} for $B > \Bchange$ are summarized in Table \ref{tab:exponents} together with the prefactors $\tC$ found in those fits.

\begin{table}[ht]
    \centering
    % 1 significant digit in error
    \begin{tabular}{c | c c c }
         $q$ & exponent in \eqref{estimate} & \multicolumn{1}{c}{measured exponent $\talpha$} & measured prefactor $\tC$ \\
         \hline
         $4$ & $12$ & $11.88           \pm 0.03$ & $2.9\cdot10^{-15}$\\
         $5$ & $10$ & $\phantom{0}9.96 \pm 0.02$ & $1.1\cdot10^{-8}\phantom{^0}$\\
         $6$ & $10$ & $\phantom{0}9.92 \pm 0.02$ & $2.8\cdot10^{-6}\phantom{^0}$\\
         $9$ & $12$ & $11.94           \pm 0.04$ & $7.0\cdot10^{-4}\phantom{^0}$
    \end{tabular}
    \hide{
        % 2 significant digits in error
        \begin{tabular}{c | c c c }
             $q$ & exponent in \eqref{estimate} & \multicolumn{1}{c}{measured exponent $\talpha$} & measured constant $\tC$ \\
             \hline
             $4$ & $12$ & $11.879 \pm 0.029$ & $2.93\cdot 10^{-15} \divideontimes 1.17$\\
             $5$ & $10$ & $\phantom{0}9.955 \pm 0.024$ & $1.08\cdot 10^{-8} \divideontimes 1.13$\\
             $6$ & $10$ & $ \phantom{0}9.916 \pm 0.025 $ & $2.80\cdot 10^{-6} \divideontimes 1.13$\\
             $9$ & $12$ & $11.940 \pm 0.042$ & $6.97\cdot 10^{-4} \divideontimes 1.22$
        \end{tabular}

        % full values
        \begin{tabular}{c | c r c }
             $q$ & exponent in \eqref{estimate} & \multicolumn{1}{c}{measured exponent $\talpha$} & measured constant $\tC$ \\
             \hline
             $4$ & $12$ & $11.87806 \pm 0.02915$ & $2.92897e-15 \divideontimes 1.16916e+00$\\
             $5$ & $10$ & $9.95533 \pm 0.02416$ & $1.08002e-08 \divideontimes 1.13378e+00$\\
             $6$ & $10$ & $ 9.91570 \pm 0.02493 $ & $2.80385e-06 \divideontimes 1.13446e+00$\\
             $9$ & $12$ & $11.93974 \pm 0.04237$ & $6.96448e-04 \divideontimes 1.22523e+00$
        \end{tabular}
    }
    \caption{Exponents in the upper bound \eqref{estimate}, exponents $\talpha$ found by fitting ansatz \eqref{fit} to the data in Figure \ref{fig:branches} for $B > \Bchange$ and the corresponding prefactors $\tC$ found in these fits for different $q$.} 
    \label{tab:exponents}
\end{table}

It is an interesting question whether there is any relation between the local maximizers $\tuE$ of Problem \ref{pb:maxdEdt} \cite{luDoering08,ayalaProtas2017} and the local maximizers $\tuB$ of Problem \ref{hilbertProblem} found here. By performing direct evaluation for different $q$ and large $\E_0$, $B$, we note that $\G(\tuB) < 0$ and $\R_q(\tuE) < 0$, which means that the states maximizing $(d/dt) \| \u \|_q^q$ can only produce $(d/dt) \| \bnabla \u \|_2^2 < 0$ and vice versa. We therefore conclude that, at least with the local maximizers of Problems \ref{pb:maxdEdt} and \ref{hilbertProblem} found in \cite{luDoering08,ayalaProtas2017} and here, it is not possible to simultaneously saturate estimates \eqref{jwekrnwer} and \eqref{estimate} with the same family of states.

%\subsubsection{structure of optimal states}
\subsection{Structure of the Optimal States}

\label{subsectionStructure}
\begin{figure}[htp]% allows for here, top or as float page, fixes long caption problem -> uses float page
\def\morphPicWidth{\dimexpr0.195\textwidth\relax}
\def\morphGapWidth{\dimexpr0.005\textwidth\relax}
    \begin{subfigure}[t]{\morphPicWidth}
        \centering
        \includegraphics[width=\textwidth]{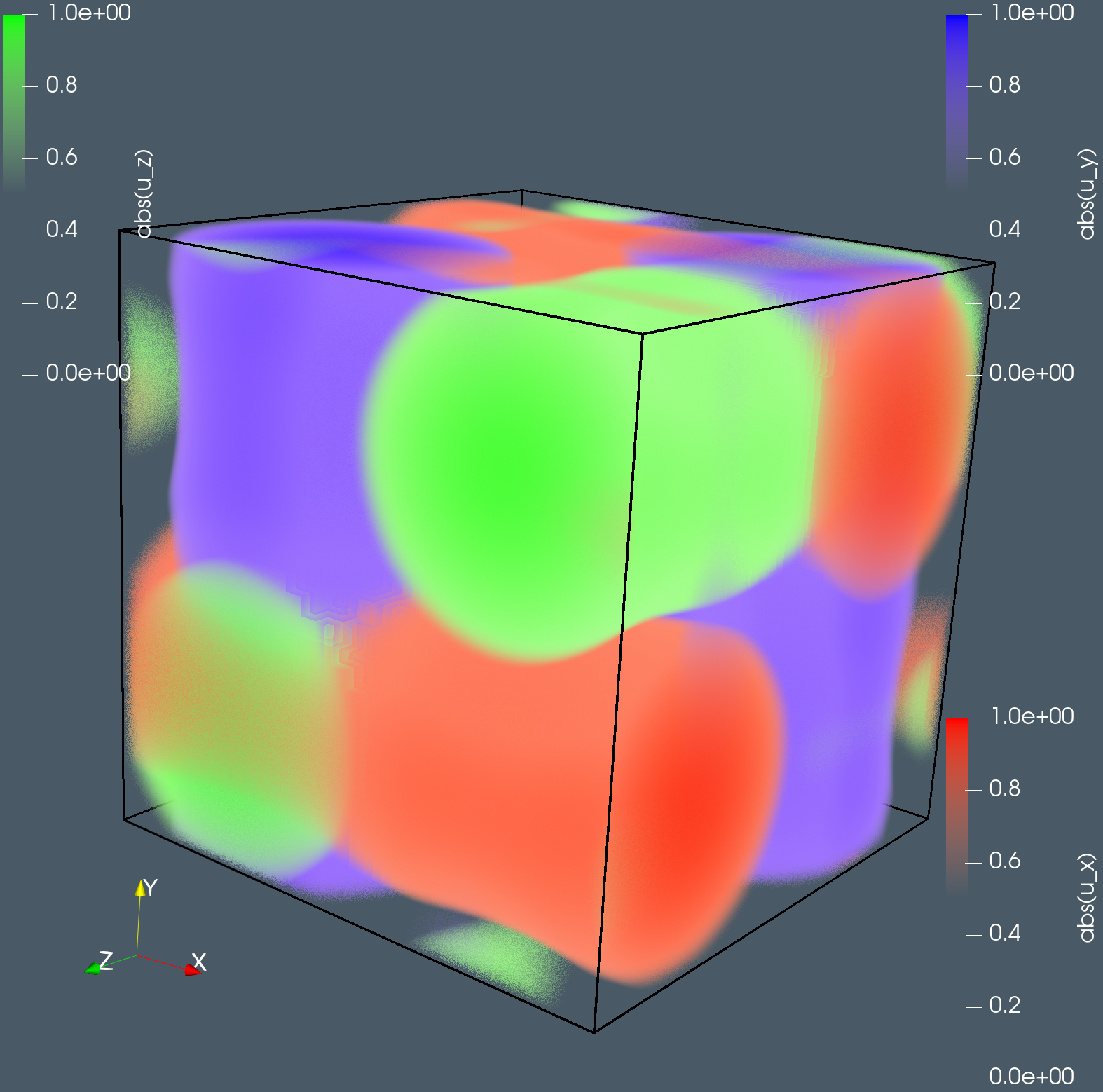}
        \caption{$B=10^0$}
        \label{fig:morphing:vel_0000}
    \end{subfigure}%
    \hspace{\morphGapWidth}%
    \begin{subfigure}[t]{\morphPicWidth}
        \centering
        \includegraphics[width=\textwidth]{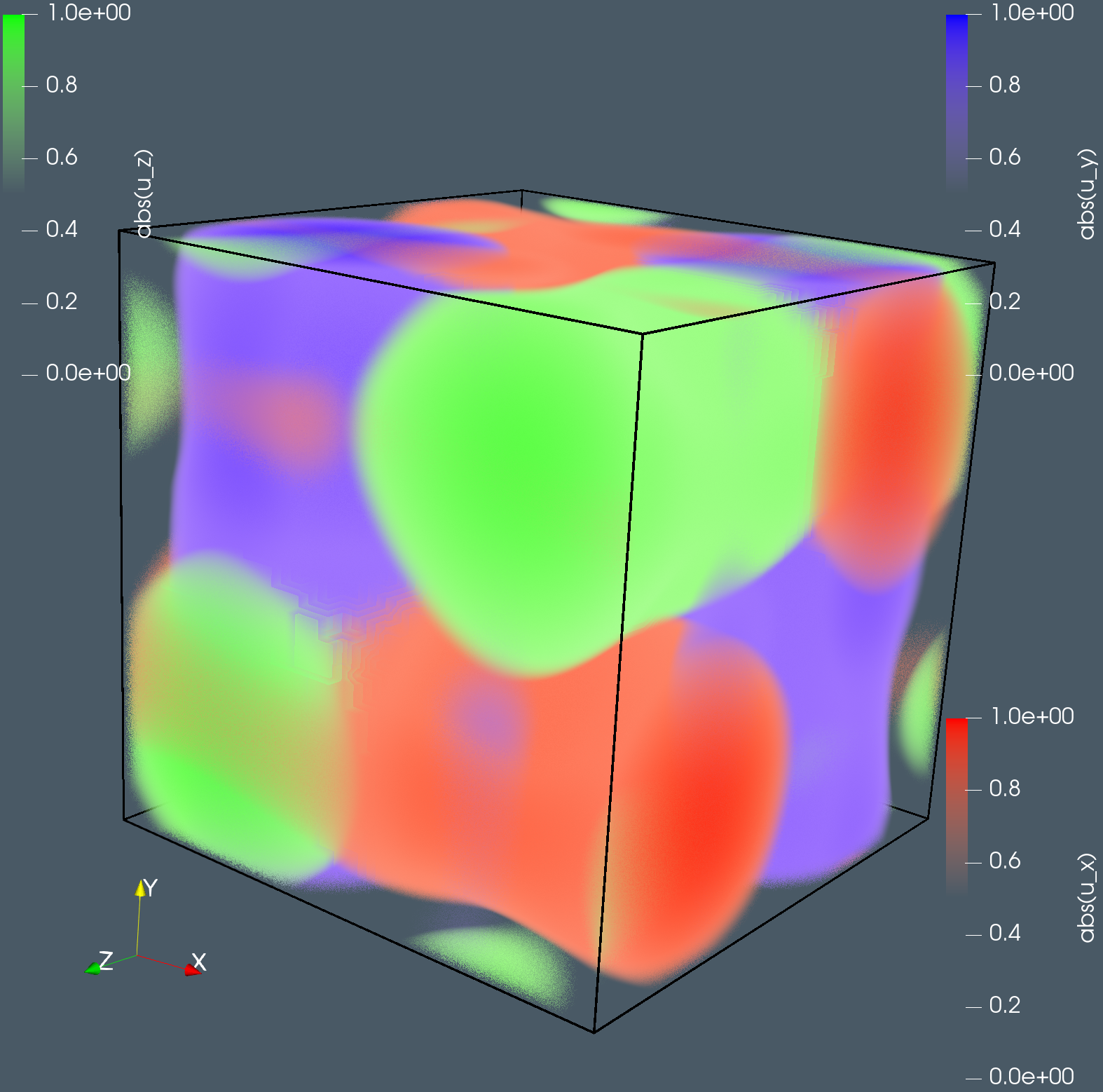}
        \caption{$B=10^{11/8}$}
        \label{fig:morphing:vel_0011}
    \end{subfigure}%
    \hspace{\morphGapWidth}%
    \begin{subfigure}[t]{\morphPicWidth}
        \centering
        \includegraphics[width=\textwidth]{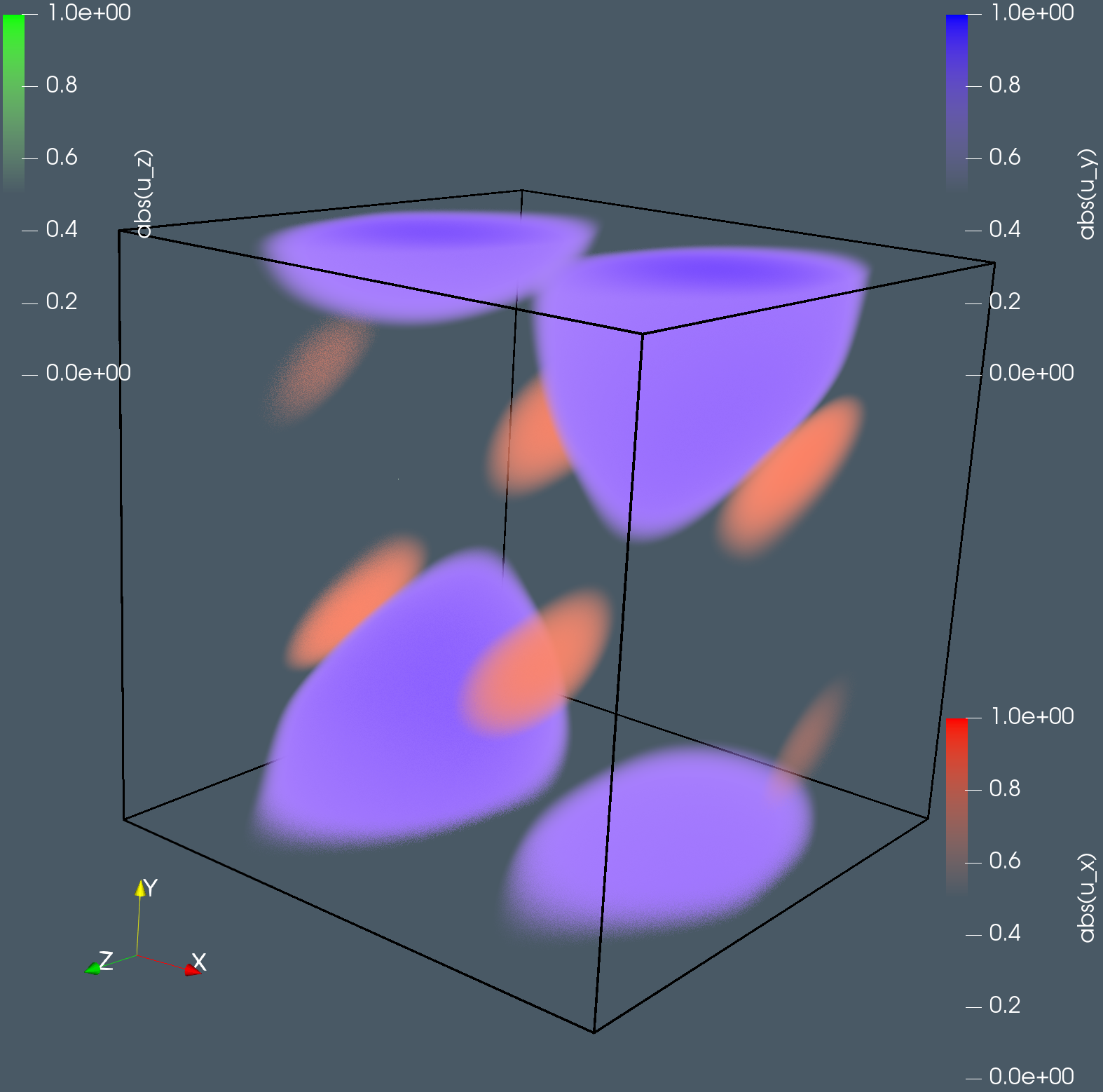}
        \caption{$B=10^{15/8}$}
        \label{fig:morphing:vel_0015}
    \end{subfigure}%
    \hspace{\morphGapWidth}%
    \begin{subfigure}[t]{\morphPicWidth}
        \centering
        \includegraphics[width=\textwidth]{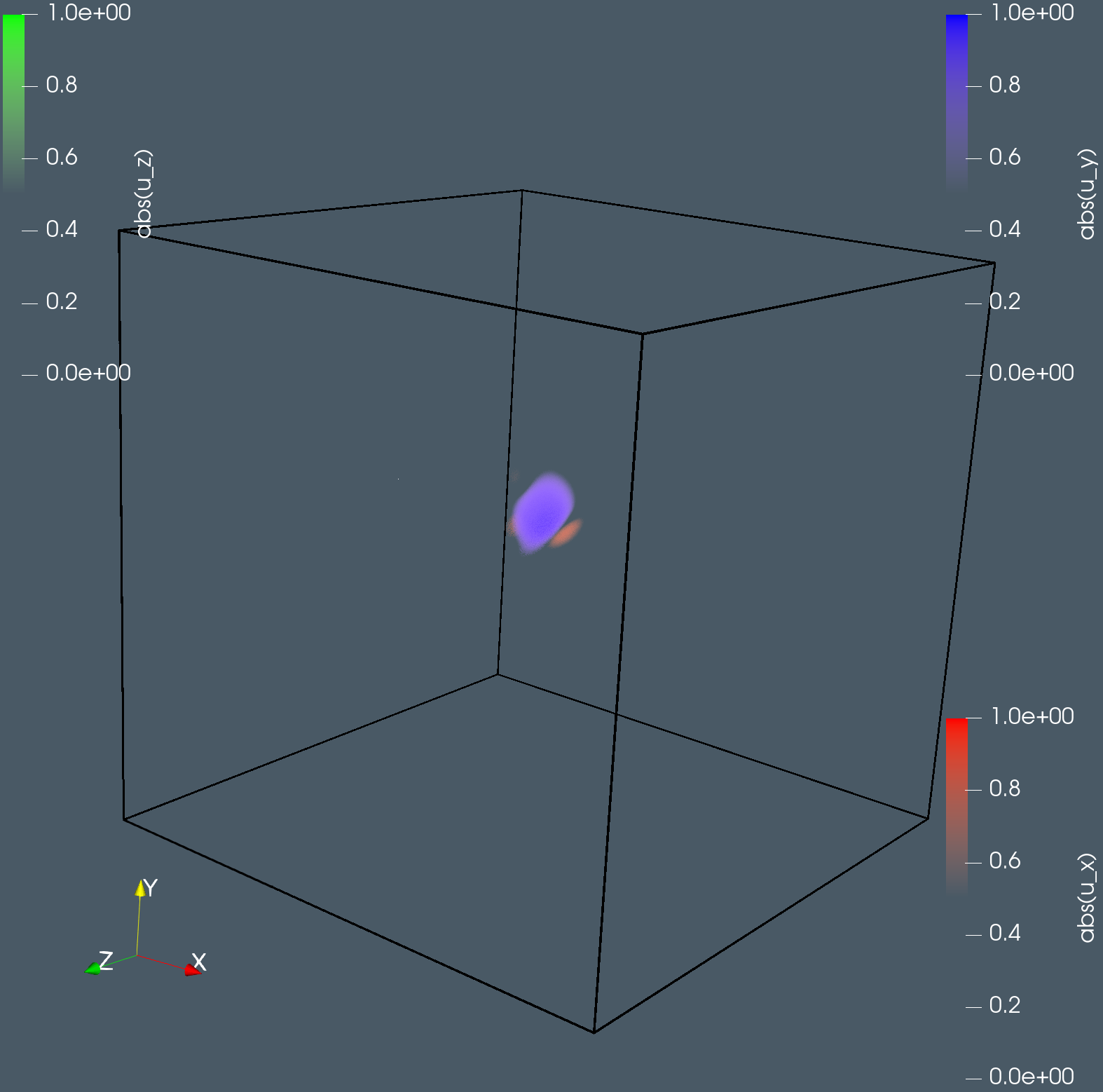}
        \caption{$B=10^{17/8}$}
        \label{fig:morphing:vel_0017}
    \end{subfigure}%
    \hspace{\morphGapWidth}%
    \begin{subfigure}[t]{\morphPicWidth}
        \centering
        \includegraphics[width=\textwidth]{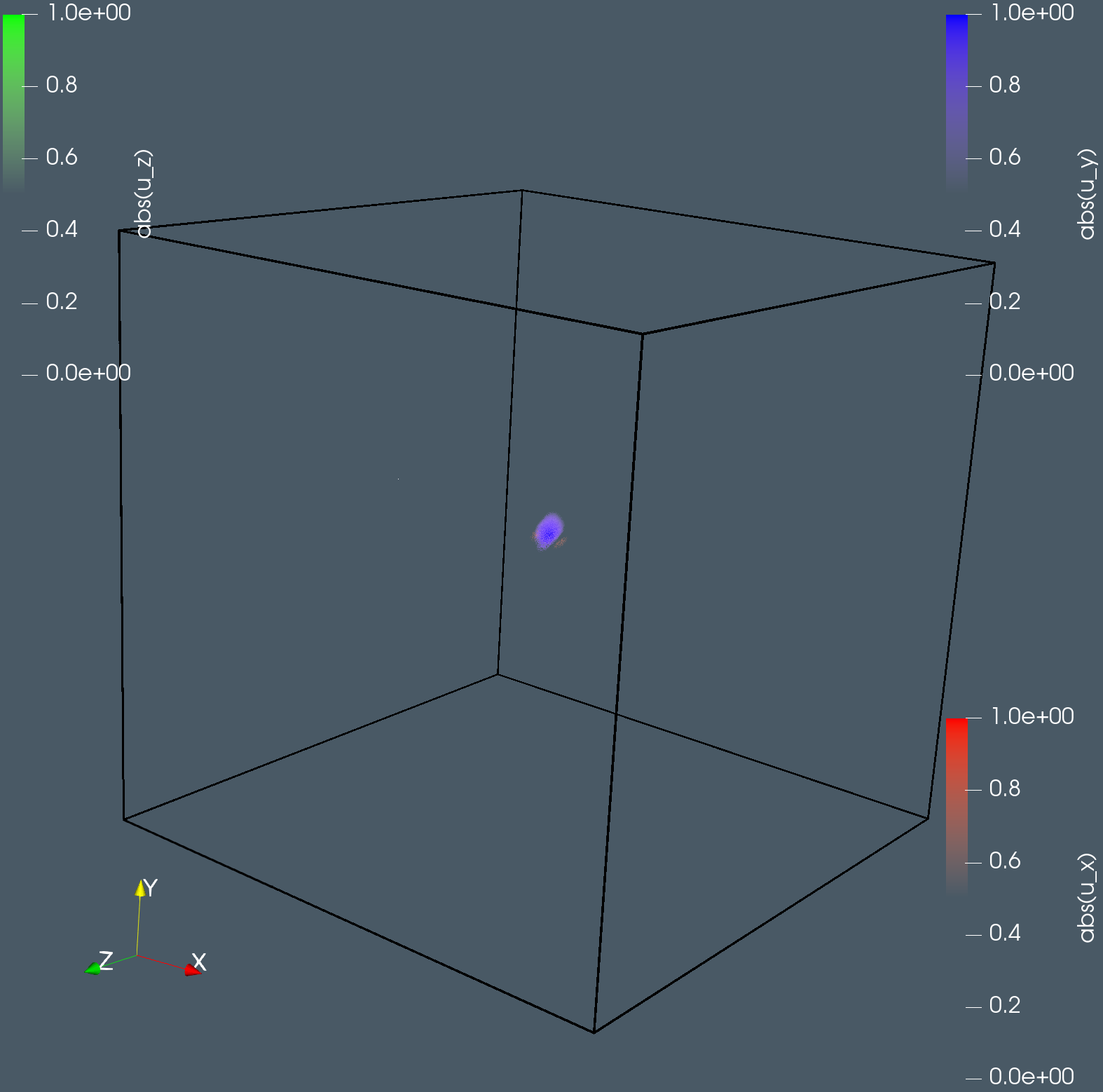}
        \caption{$B=10^{18/8}$}
        \label{fig:morphing:vel_0018}
    \end{subfigure}%
    \vspace{-0.5em}
    \begin{center}
    %    \small{Velocity}
    \end{center}    
    \begin{subfigure}[t]{\morphPicWidth}%
        \centering
        \includegraphics[width=\textwidth]{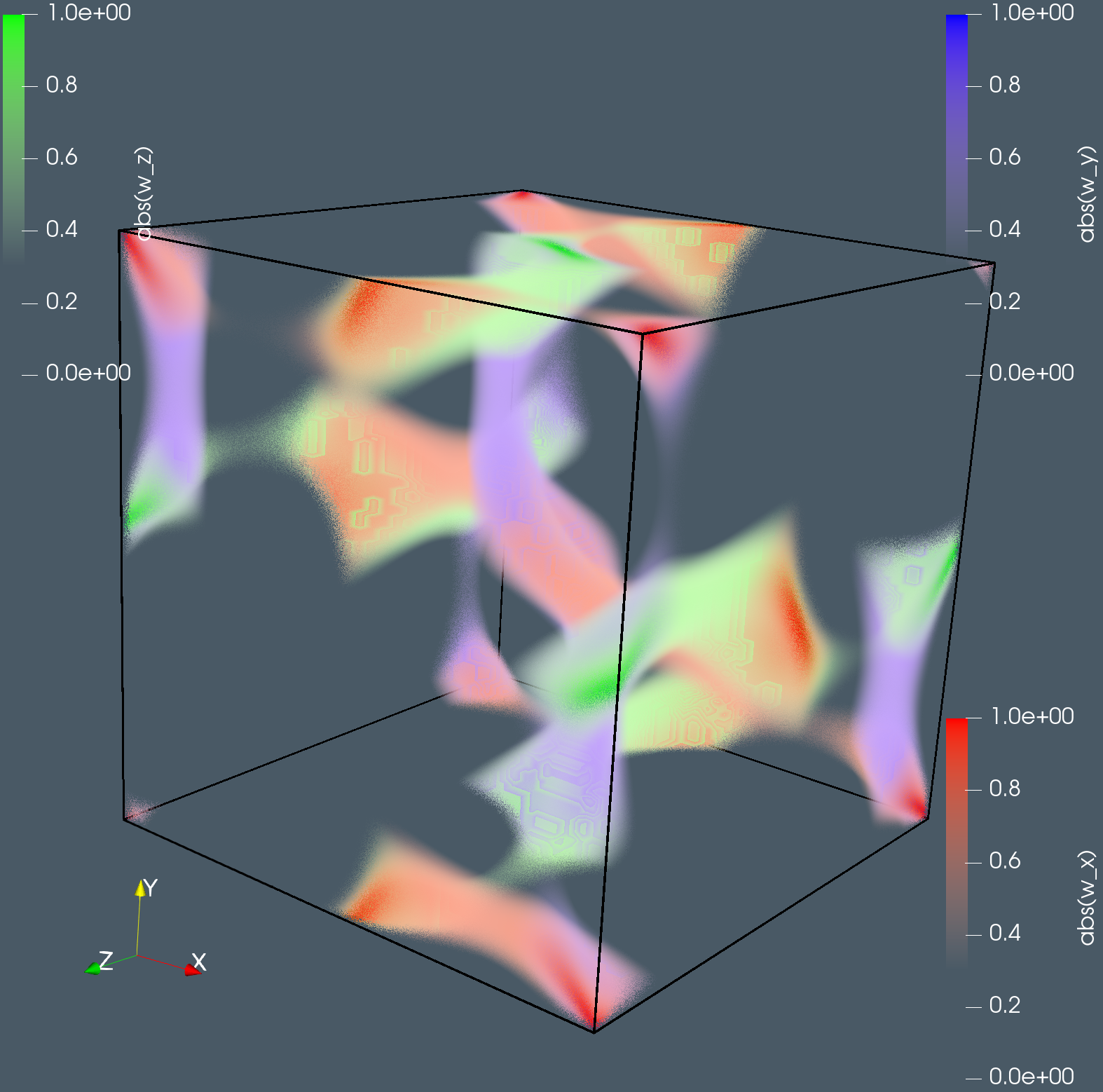}
        \caption{$B=10^0$}
        \label{fig:morphing:vort_0000}
    \end{subfigure}%
    \hspace{\morphGapWidth}%
    \begin{subfigure}[t]{\morphPicWidth}%
        \centering
        \includegraphics[width=\textwidth]{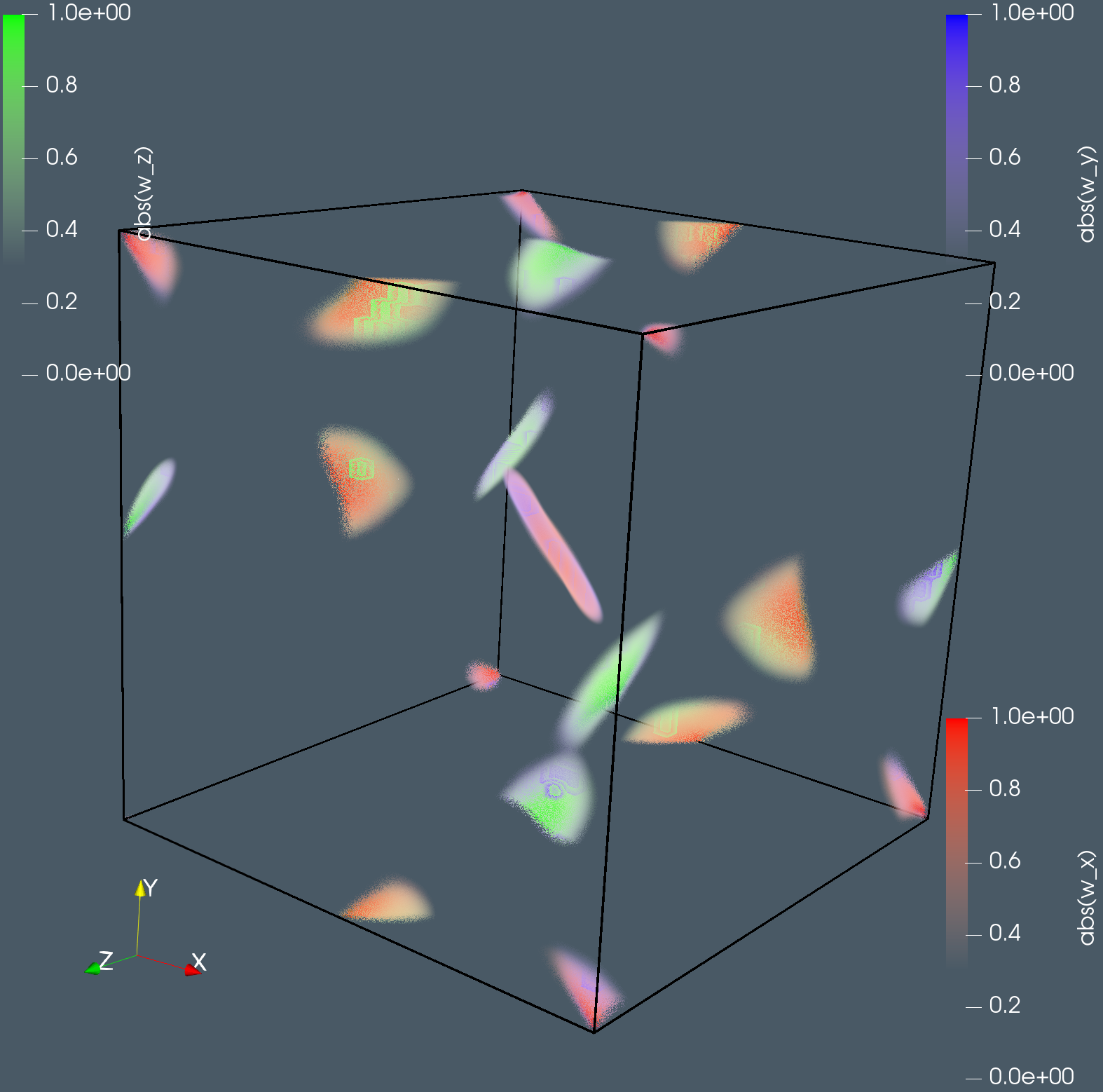}
        \caption{$B=10^{11/8}$}
        \label{fig:morphing:vort_0011}
    \end{subfigure}%
    \hspace{\morphGapWidth}%
    \begin{subfigure}[t]{\morphPicWidth}%
        \centering
        \includegraphics[width=\textwidth]{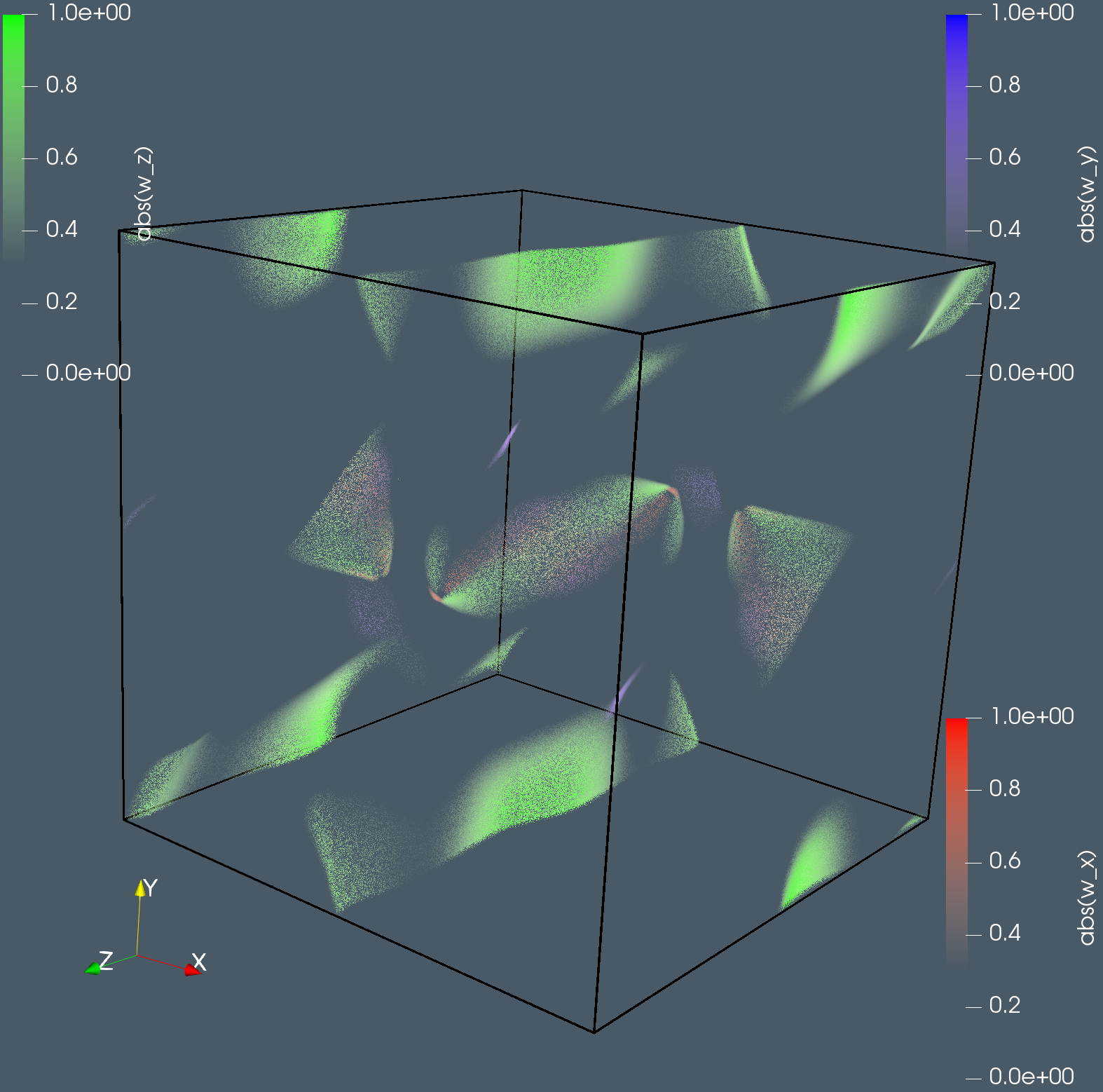}
        \caption{$B=10^{15/8}$}
        \label{fig:morphing:vort_0015}
    \end{subfigure}%
    \hspace{\morphGapWidth}%
    \begin{subfigure}[t]{\morphPicWidth}%
        \centering
        \includegraphics[width=\textwidth]{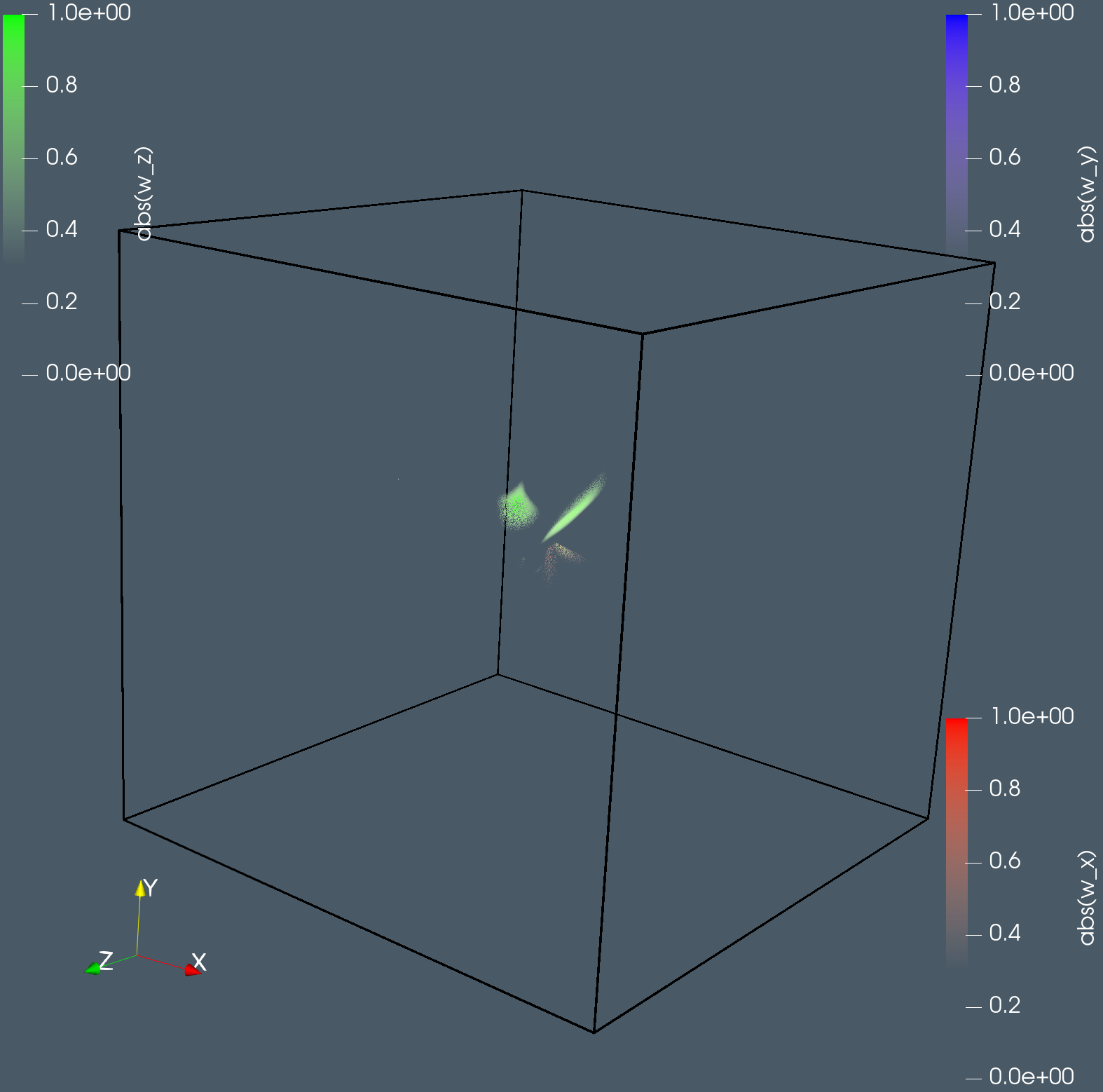}
        \caption{$B=10^{17/8}$}
        \label{fig:morphing:vort_0017}
    \end{subfigure}%
    \hspace{\morphGapWidth}%
    \begin{subfigure}[t]{\morphPicWidth}%
        \centering
    \includegraphics[width=\textwidth]{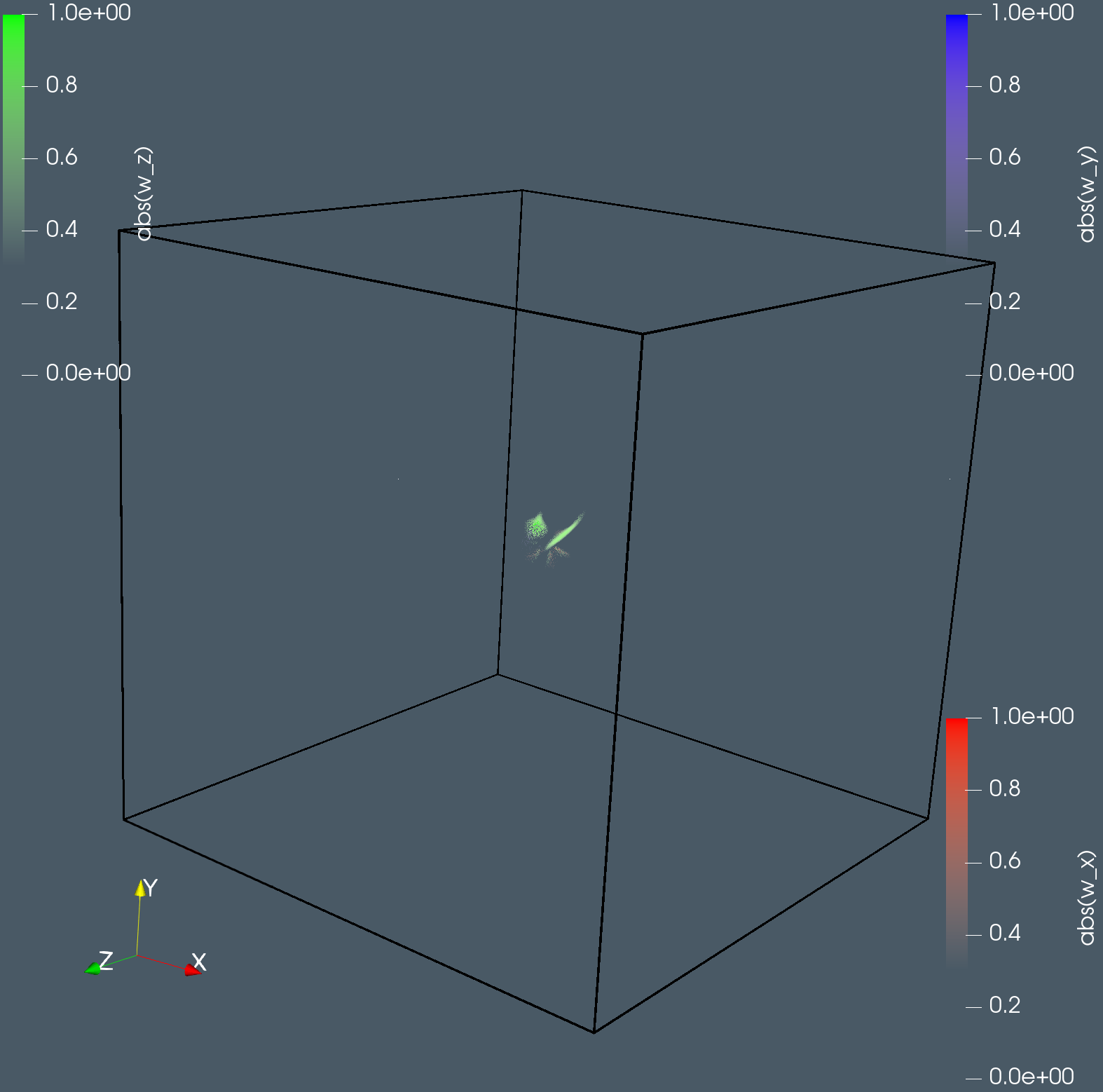}
        \caption{$B=10^{18/8}$}
        \label{fig:morphing:vort_0018}
    \end{subfigure}%
    \vspace{-0.5em}
    \begin{center}
    %    \small{Vorticity}
    \end{center}
    \vspace{-1em}
        \caption{Normalized  (\subref{fig:morphing:vel_0000})--(\subref{fig:morphing:vel_0018}) velocity $\tuB(\x)$ and (\subref{fig:morphing:vort_0000})--(\subref{fig:morphing:vort_0018}) vorticity $(\bnabla \times \tuB)(\x)$ fields in the solutions of Problem \ref{hilbertProblem} with $q=5$ and indicated values of the  constraint parameter. The different colors represent the Cartesian components of the fields:  {\color{paraviewRed}(red)} $x_1$, {\color{paraviewBlue}(blue)} $x_2$ and {\color{paraviewGreen}(green)} $x_3$. Panels (\subref{fig:morphing:vel_0000}), (\subref{fig:morphing:vel_0011}), (\subref{fig:morphing:vort_0000}), and (\subref{fig:morphing:vort_0011}) show the \smallB{} regime, panels (\subref{fig:morphing:vel_0015}), (\subref{fig:morphing:vort_0015}) correspond to the transition phase whereas panels (\subref{fig:morphing:vel_0017}),  (\subref{fig:morphing:vel_0018}),  (\subref{fig:morphing:vort_0017}), and (\subref{fig:morphing:vort_0018}) represent the \largeB{} regime.}
    \label{fig:morphing}
\end{figure}

The physical-space structure of the optimal fields $\tuB$ found by solving Problem \ref{hilbertProblem} changes as the constraint parameter $B$ increases and the key features of these flow patterns are shared by the solutions obtained for different values of $q$. As $B$ increases from zero, the maximizers gradually deviate from the Arnold-Beltrami-Childress flow \eqref{ABC}, which is the solution of Problem \ref{hilbertProblem} for $B \rightarrow 0$, cf.~Section \ref{sectionSmallDataLimit}. For intermediate values of $B$ the structure of $\tuB$ undergoes a transition and the flow that emerges as a result persists for increasing values of $B$ while shrinking in size such that it becomes more localized. This evolution of $\tuB$ as $B$ increases from "small" to "large" values is illustrated for $q = 5$ in Figure \ref{fig:morphing}.

\begin{figure}[htp]% allows for here, top or as float page, fixes long caption problem -> uses float page
\centering
    \begin{subfigure}[t]{0.35\textwidth}
        \centering
        \includegraphics[width=\textwidth]{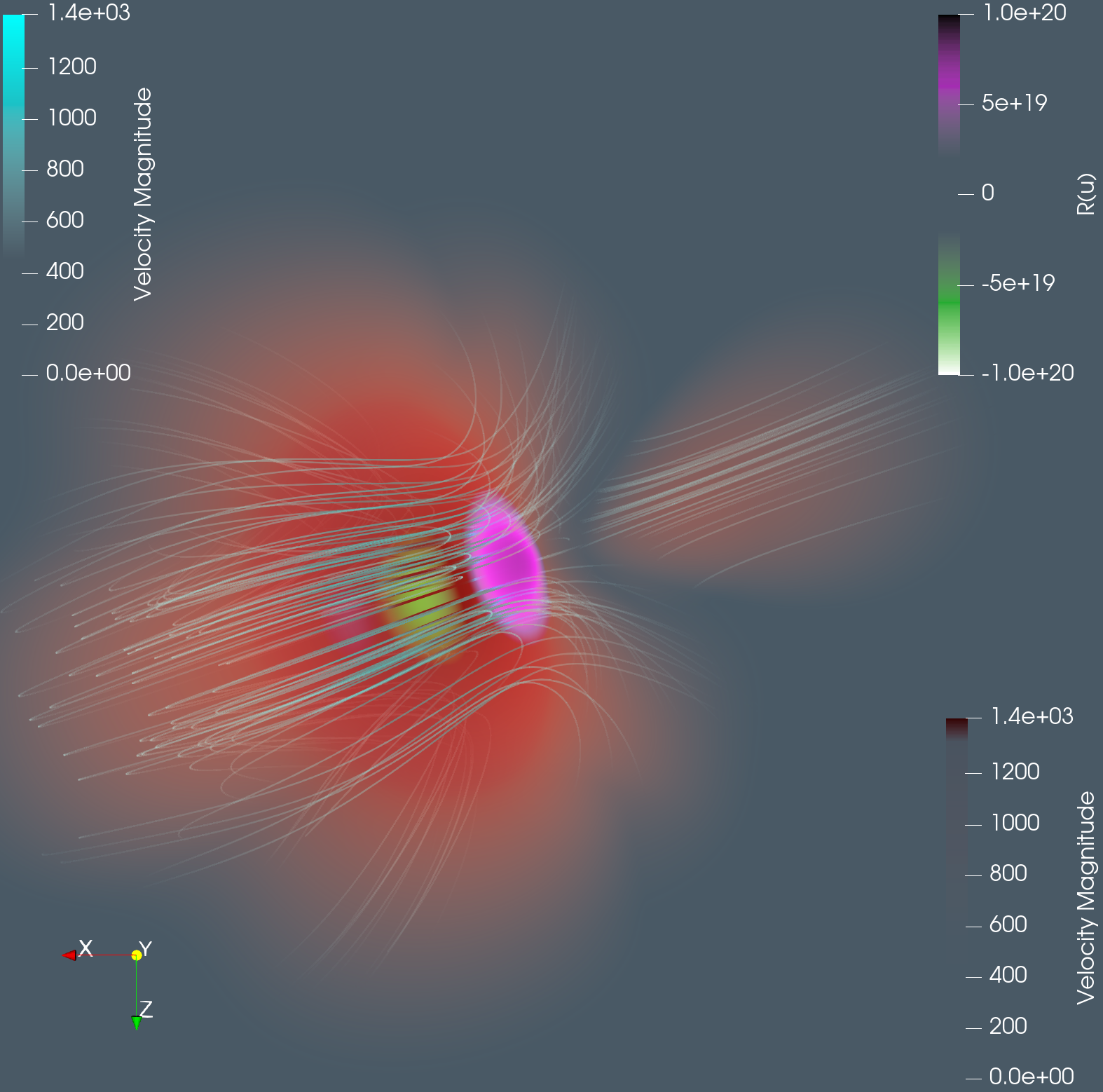}
        \caption{}
        %\caption{{\color{paraviewRed}Velocity magnitude}, {\color{paraviewGreen}stream lines}, and pointwise objective functional}
        \label{fig:structure_q5_B19:vel}
    \end{subfigure}%
    \hspace{0.1\textwidth}%
    \begin{subfigure}[t]{0.35\textwidth}%
        \centering
        \includegraphics[width=\textwidth]{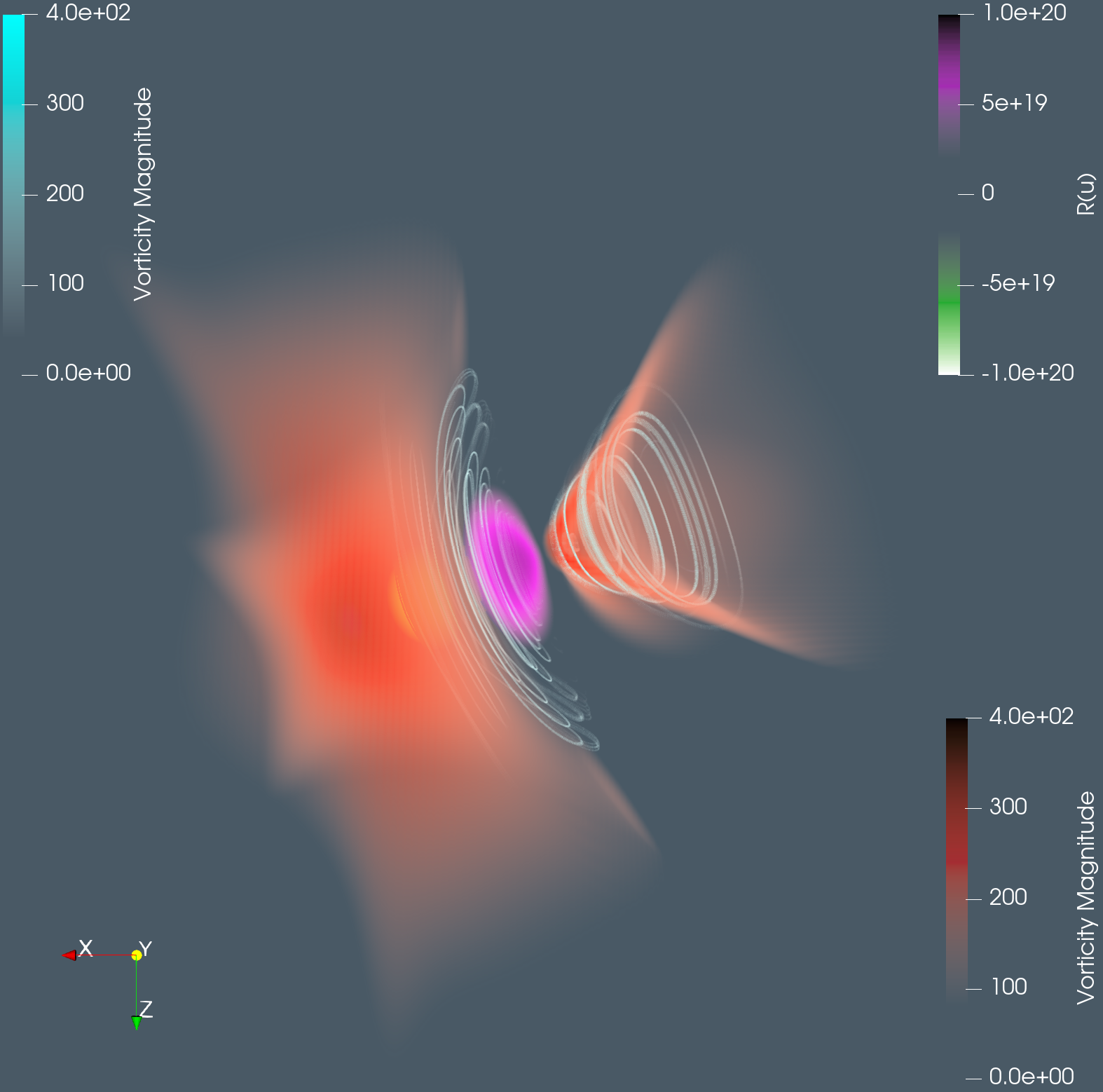}
        \caption{}
        %\caption{{\color{paraviewRed}Vorticity magnitude}, {\color{paraviewGreen}vortex lines}, and pointwise objective functional}
        \label{fig:structure_q5_B19:vort}
    \end{subfigure}%
%    \hspace{0.04\textwidth}%
%    \begin{subfigure}[t]{0.30\textwidth}%
%        \centering
%        \includegraphics[width=\textwidth]{pictures/structure_of_fields/q5_B19_qCrit_Ru.png}
%        \caption{}
%        \label{fig:structure_q5_B19:qCrit}
%    \end{subfigure}%
    \caption{Structure of the solution $\tuB$ of the Problem \ref{hilbertProblem} with $q=5$ and $B = 10^{9/4}\approx 177.8$: (\subref{fig:structure_q5_B19:vel}) {\color{paraviewRed}red} represents the velocity magnitude $|\tuB(\x)|$ and {\color{paraviewTurquoise}turquoise} the stream lines, (\subref{fig:structure_q5_B19:vort}) {\color{paraviewRed}red} represents the vorticity magnitude $|\bnabla \times \tuB(\x)|$ and {\color{paraviewTurquoise}turquoise} the vortex lines  (i.e., lines tangent to the vorticiy field).  %{\color{paraviewBlue}Positive} and {\color{paraviewOrange}negative} values of the Q-criterion are shown in Figure (\subref{fig:structure_q5_B19:qCrit}). 
    Positive and negative values of the integrand expression \eqref{r} are marked with {\color{paraviewPurple}purple} and {\color{paraviewGreen}green} in both panels.}
    \label{fig:structure_q5_B19}
\end{figure}

We now take a closer look at the maximizers found for large $B$ and in Figure \ref{fig:structure_q5_B19} we show $\tuB$ obtained for $q = 5$ and $B = 10^{9/4} \approx 177.8$, which is the representative case already discussed in Subsection \ref{subsectionRepresentativeValues}. In the figure we also plot the integrand expression of the objective functional \eqref{defRq} evaluated at $\tuB$, i.e.,
\begin{align}
    \integrandFull(\tuB) &:= - \nu |\tuB|^{q-2}|\bnabla \tuB|^2 - \frac{4(q-2)\nu}{q^2} \Big|\bnabla |\tuB|^{\frac{q}{2}}\Big|^2 \nonumber
    \\
    &\qquad\qquad\qquad - (q-2) |\tuB|^{q-4} \tuB\cdot (\tuB\cdot \bnabla) \tuB\, \Delta^{-1}\left[\bnabla \tuB\mathcolon \bnabla (\tuB)^T\right]
\label{r}
\end{align}
such that $\R_q(\tuB) = \int_\Omega \integrandFull(\tuB(\x))\ d\x$. We see that, in contrast to solutions of Problem \ref{pb:maxdEdt} which have a clear form of two colliding nearly axisymmetric vortex rings, the structure of our maximizers $\tuB$ is more complex and cannot be reduced to  simple flows. However, they also appear to be nearly axisymmetric although without a symmetry plane. It is interesting to note that the integrand expression $r(\tuB(\x))$ attains its maximum value at the center of the structure. We add that, for large $B$, the physical-space structure of the maximizers $\tuB(\x)$ becomes more complex for increasing values of $q$.

%\subsubsection{\texorpdfstring{$q\to 3, q=3$}{q to 3 and q equal 3}}

\subsection{The Limit \texorpdfstring{$q\to 3$}{q to 3}}
\label{subsectionLimitQto3}
As discussed in Section \ref{sectionIntro}, the case with $q=3$ is potentially the most interesting due to the scaling invariance of the Navier-Stokes equations in $\RR^3$ and the breakdown of the bound \eqref{estimate}. In order to shed light on this issue, we have attempted to solve Problem \ref{hilbertProblem} for $q = 3$. However, while solutions $\tuB$  such that $\R_3(\tuB) < 0$ could be found in the small-data regime, cf.~Figure \ref{fig:branches}(b), the iterations \eqref{RCG1}--\eqref{RCG2} would not converge for $B > \Bchange$. More specifically, the objective functional $\R_3(\u_k)$ grew without bound as $k \rightarrow \infty$ with the corresponding approximations $\u_k$ requiring increasing resolutions to be accurately represented, which is a manifestation of the ill-posedness of the optimization problem. This behavior gradually emerged when Problem \ref{hilbertProblem} was solved for $q \rightarrow 3$ with a modification of Algorithm \ref{alg1} where $B$ is kept constant and $q$ decremented from 4. Performing such calculations for two fixed values of the constraint parameter $B$ made it possible to determine the exponents $\talpha$ in ansatz \eqref{fit}. They are reported in Table \ref{tab:exponents_qTo3} indicating divergence of  the upper bound in \eqref{estimate} as $q$ is reduced down from 4. However, these results should be interpreted with caution due to relatively large uncertainties. 

 %This indicates that indeed the left-hand side of the a priori bound diverges. However, these simulations showed a similar non-convergence as the case $q=3$, although for slightly higher constraint values such that we have $2$ data points for each value of $q$ allowing us to just barely calculate the scaling exponents. These exponents are given in Table \ref{tab:exponents_qTo3}, but interpreted with caution due to relatively large uncertainties. 

\begin{table}[ht]
    \centering
    \begin{tabular}{c | c c}
         $q$ & exponent in \eqref{estimate} & measured exponent $\talpha$ \\
         \hline
         %3.938114489599 & 12.333922321716292 & 12.26424962126599\\
         %3.877186433297 & 12.71723750140159 & 12.499938704999725\\
         %3.81720101796 & 13.159335778117244 & 12.834695416672794\\
         %3.758143659635 & 13.67221103162705 & 13.373478180918541\\
         %3.7 & 14.27142857142857 & 14.067846200233713 
         3.94 & 12.33 & 12.26\\
         3.88 & 12.72 & 12.50\\
         3.82 & 13.16 & 12.83\\
         3.76 & 13.67 & 13.37\\
         3.70 & 14.27 & 14.07 
    \end{tabular}
    \caption{Exponents in the upper bound \eqref{estimate} and the exponents $\talpha$ found by fitting ansatz \eqref{fit} to the data obtained by solving Problem \ref{hilbertProblem} with fixed values of $B$ and $q$ decremented from 4.}
    \label{tab:exponents_qTo3}
\end{table}

\section{Conclusions}
\label{sectionConclusion}

%\subsubsection{the bounds seem to be sharp}

Motivated by the Ladyzhenskaya-Prodi-Serrin conditional regularity results \eqref{lps}, we have investigated the sharpness of estimate \eqref{estimate} on the rate of growth of the $L^q$ norm of the velocity field. This is a companion problem with respect to the question about the sharpness of estimate \eqref{jwekrnwer} on the rate of growth of enstrophy which is related to the enstrophy condition \eqref{maxEinf} and was studied ealier in \cite{luDoering08,ayalaProtas2017}, cf.~Problem \ref{pb:maxdEdt}.
In particular, we formulated a family of optimization problems where the instantaneous rate of growth of the $L^q$ norm is maximized with respect to the velocity field subject to suitable constraints. They have two closely-related versions, Problems \ref{banachProblem} and  \ref{hilbertProblem}, stated in Sobolev spaces without and with  
the Hilbert structure.  However, only the latter formulation is amenable to numerical solution  and was considered in this study. Given its Riemannian structure, it was solved using a state-of-the-art Riemannian conjugate-gradient method in the optimize-then-discretize setting where the gradient direction was computed using methods of calculus of variations.

Branches of local maximizers of Problem \ref{hilbertProblem} were found for $q = 4,5,6,9$. As the main finding of our study, the rate of growth $(d/dt) \| \u(t) \|_q^q$
evaluated on each of these branches saturates the upper bound in estimate \eqref{estimate} as $B \rightarrow \infty$. This demonstrates that this bound is sharp and therefore cannot be fundamentally improved, except for perhaps a refinement of the constant prefactor. This is a nontrivial conclusion since while it is known that the estimates used at the different steps in the derivation of \eqref{estimate}, cf.~\ref{appendixProofOfAprioriBound}, are sharp, it is not a priori obvious this property should be preserved when these estimates are chained together and evaluated on a single field. Since estimate \eqref{estimate} was already saturated by solutions of Problem \ref{hilbertProblem}, there was no need to solve the more general Problem \ref{banachProblem} which would have been computationally intractable.

Unlike the local maximizers $\tuE$ of Problem \ref{pb:maxdEdt}, which have the clear form of two colliding axisymmetric vortex rings \cite{luDoering08,ayalaProtas2017}, the maximizers $\tuB$ of Problem \ref{hilbertProblem} found for different $q$ and $B > \Bchange$ do not reveal a simple structure in the physical space, cf.~Section \ref{subsectionStructure}. However, the two families of maximizers share the property that they become increasingly localized as $B \rightarrow \infty$. Interestingly, we have $\G(\tuB) < 0$ and $\R_q(\tuE) < 0$, indicating that the states saturating one of the bounds are in fact quite far from simultaneously also saturating the other bound. Since at a hypothetical blow-up both, \eqref{lps} and \eqref{maxEinf} need to be violated, this suggests that none of these maximizers could give rise to singularity formation in finite time when used as the initial condition $\u_0$ for the Navier-Stokes system \eqref{nse}--\eqref{incompressible}.

Finally, when $q = 3$, the inapplicability of bound \eqref{estimate} together with the observation that Problem \ref{hilbertProblem} may not be well posed in this case suggest that there is likely no finite a priori bound on $(d/dt) \| \u(t) \|_3$. This is a surprising conclusion since in various extreme flows $\| \u(t) \|_3$ has revealed a weaker growth than $\| \u(t) \|_q$, $q > 3$ \cite{Hou:22:NS, ramirez2025}.

%\subsubsection{\texorpdfstring{$\lim_{q\to 3}$}{q to 3} does not seem to exist}

\newpage
\appendix

\section{Proof of the A Priori Bound\texorpdfstring{ \eqref{estimate}}{}}
\label{appendixProofOfAprioriBound}

Here we follow the ideas of \cite{robinsonSadowski2014}. Using \eqref{p} and applying Hölder's and Young's inequalities to \eqref{defRq} yields
\begin{align}
    \frac{1}{q}\frac{d}{dt}\|\u\|_q^q + \nu & \int_{\Omega} |\u|^{q-2}|\bnabla \u|^2 \, d\x 
    \\
    &\leq (q-2)\int_\Omega p |\u|^{q-4} \u\cdot (\u\cdot \bnabla) \u \, d\x
    \\
    &\leq (q-2) \left\| p |\u|^{\frac{q-2}{2}}\right\|_2 \left\| |\u|^{\frac{q-2}{2}} |\bnabla \u| \right\|_2
    \\
    &\leq \frac{(q-2)^2}{2\nu} \left\| p |\u|^{\frac{q-2}{2}}\right\|_2^2 + \frac{\nu}{2}\int_{\Omega} |\u|^{q-2} |\bnabla \u|^2 \, d\x,
\end{align}
which implies
\begin{align}
    \frac{1}{q}\frac{d}{dt}\|\u\|_q^q + \frac{\nu}{2}  \int_{\Omega} |\u|^{q-2}|\bnabla \u|^2 \, d\x
    &\leq \frac{(q-2)^2}{2\nu} \int_{\Omega} |p|^2 |\u|^{q-2} d\x.\label{kjwne}
\end{align}
To further estimate the right-hand side note that the pressure satisfies (see \cite[Lemma 3]{robinsonSadowski2014})
\begin{align}
    \|p\|_{r} \leq C \|\u\|_{2r}^2.\label{ajdfnkkjsf}
\end{align}
%For a proof of \eqref{ajdfnkkjsf} see robinsonRodrigoSadowski Lemma 5.1 + proof in Appendix B or \cite[Lemma 3]{robinsonSadowski2014}. 
Therefore, Hölder's inequality, \eqref{kjwne} and \eqref{ajdfnkkjsf} yield
\begin{align}
    \frac{1}{q}\frac{d}{dt}\|\u\|_q^q + \frac{\nu}{2} \int_{\Omega} |\u|^{q-2}|\bnabla \u|^2 \, d\x
    &\leq \frac{(q-2)^2}{2\nu} \int_{\Omega} |p|^2 |\u|^{q-2} d\x 
    \\
    &\leq \frac{(q-2)^2}{2\nu} \|p^2\|_{\frac{q+2}{4}} \left\| |\u|^{q-2} \right\|_{\frac{q+2}{q-2}}
    %\\
    %&= \frac{(q-2)^2}{2\nu} \|p\|_{L^\frac{q+2}{2}}^2 \| \u \|_{q+2}^{q-2}
    \\
    &\leq C\nu^{-1} \|\u\|_{q+2}^{q+2}.
    \label{njwerjkw}
\end{align}
The Gagliardo-Nirenberg type inequality (see \cite[Lemma 2]{robinsonSadowski2014})
\begin{align}
    \|\u\|_{3q}^q \leq C \int_{\Omega} |\u|^{q-2}|\bnabla \u|^2\, d\x
\end{align}
holds for all $q>2$ if $\u\in W^{1,\frac{3q}{q+1}}(\Omega)$. Therefore, by \eqref{njwerjkw}, Hölder's and Young's inequalities we have
\begin{align*}
    \frac{1}{q}\frac{d}{dt}\|\u\|_q^q + \frac{\nu}{2}  \int_{\Omega} |\u|^{q-2}|\bnabla \u|^2 \, d\x
    &\leq C\nu^{-1} \int_{\Omega} |\u|^{q+2}\, d\x
    \\
    &\leq C\nu^{-1} \big\||\u|^{q-1}\big\|_{\frac{q}{q-1}}\big\||\u|^{3}\big\|_{q}
    \\
    &= C\nu^{-1} \|\u\|_{q}^{q-1}\|\u\|_{3q}^3
    \\
    &\leq \frac{\nu}{2\widetilde{C}} \|\u\|_{3q}^{q} + C\nu^{-\frac{q+3}{q-3}} \|\u\|_q^{\frac{q(q-1)}{q-3}}
    \\
    &\leq \frac{\nu}{2}\int_{\Omega} |\u|^{q-2}|\bnabla \u|^2\, d\x + C\nu^{-\frac{q+3}{q-3}} \|\u\|_q^{\frac{q(q-1)}{q-3}}
\end{align*}
for $q>3$ and some $\widetilde{C} >0$, which implies \eqref{estimate}.

\section{Time Until Blow-up}
\label{appendixSaturationImpliesBlowup}
In this appendix we introduce a simple lemma establishing an explicit dependence of the hypothetical blow-up time on the rate of growth of $\| \u(t) \|_q$ for $q > 3$.
\begin{lemma}
    \label{lemmaSaturationImpliesBlowUp}
    Let $3<q<\infty$, $1<\alpha\leq \frac{q-1}{q-3}$, $\bar C>0$, $\|\uInitPDE\|_q>1$ and
    \begin{align}
        T^\star = \frac{1}{(\alpha-1) \bar C \|\uInitPDE\|_q^{q(\alpha-1)}}.
    \end{align}
    If $\u(t)$ satisfies $\u(0)=\uInitPDE$ and
    \begin{align}
        \frac{d}{dt} \|\u(t)\|_q^q =  \bar C \|\u(t)\|_q^{q\alpha}
        %\left(\|\u\|_q^{q}\right)^\alpha
        \label{duqdt}
    \end{align}
    for all $0\leq t<T^\star$, then
    \begin{align}
        \lim_{T\to T^\star} \int_0^T \|\u(t)\|_q^p \, dt = \infty,
    \end{align}
    where $\frac{2}{p}+\frac{3}{q}=1$.
\end{lemma}
The proof follows immediately upon integrating \eqref{duqdt} in time and makes no reference to the Navier-Stokes system \eqref{nse}--\eqref{incompressible}. For actual Navier-Stokes flows blow-up may occur only if the exponent $\alpha$ is sufficiently large, namely,
\begin{alignat}{3}
    \frac{3q-2}{3(q-2)}&< \alpha &&\leq \frac{q-1}{q-3} &\qquad&\text{for }2\leq q\leq 6
    \\
    \frac{q-2}{q-3} &< \alpha &&\leq \frac{q-1}{q-3}&\qquad&\text{for }q>6
    %\left.\begin{matrix}\text{if }2\leq q\leq 6 &\frac{3q-2}{3(q-2)}\\\text{if }6<q\phantom{<\infty}&\frac{q-2}{q-3}\end{matrix}\right\}< \alpha \leq \frac{q-1}{q-3}.
\end{alignat}
The lower bounds on $\alpha$ are established in \cite[Section 3.4]{ramirezProtas2026} and the upper bound is from \eqref{estimate}.

\section{Computation of the Gateaux Differential\texorpdfstring{ \eqref{dR}}{}}%
\label{appendixCalcRderivative}%
Using the definition of the Gateaux differential, we obtain
\begin{align}
    d&\R_q(\v;\v') \nonumber
    \\
    &= \lim_{\epsilon\to 0}\frac{1}{\epsilon}\left[\R_q(\v+\epsilon \v')-\R_q(\v)\right]
    \\
    %&= \lim_{\epsilon\to 0} \frac{1}{\epsilon}\Bigg(
    %\int_\Omega |(\v+\epsilon \v')|^{q-2} (\v+\epsilon \v')\cdot \bnabla \Delta^{-1} \left(\bnabla (\v+\epsilon \v')\mathcolon\bnabla (\v+\epsilon \v')^T\right) d\x
    %\\
    %&\qquad\qquad\quad + \nu \int_{\Omega} |(\v+\epsilon \v')|^{q-2}(\v+\epsilon \v')\cdot \Delta (\v+\epsilon \v') d\x
    %\\
    %&\qquad\qquad\quad - 
    %\int_\Omega |\v|^{q-2} \v\cdot \bnabla \Delta^{-1} \left(\bnabla \v\mathcolon\bnabla \v^T\right) \, d\x - \nu \int_{\Omega} |\v|^{q-2}\v\cdot \Delta \v \, d\x\Bigg)
    %\\
    &= (q-2) \int_\Omega |\v|^{q-4} (\v\cdot \v') \v \cdot \bnabla \Delta^{-1} \left(\bnabla \v\mathcolon\bnabla \v^T\right) d\x \nonumber
    \\
    &\qquad + \int_\Omega |\v|^{q-2} \v' \cdot \bnabla \Delta^{-1} \left(\bnabla \v\mathcolon\bnabla \v^T\right) d\x \nonumber
    \\
    &\qquad + 2 \int_\Omega |\v|^{q-2} \v \cdot \bnabla \Delta^{-1} \left(\bnabla \v'\mathcolon\bnabla \v^T\right) d\x \nonumber
    \\
    &\qquad + (q-2) \nu \int_{\Omega} |\v|^{q-4} (\v\cdot\v')\v\cdot \Delta \v \, d\x + \nu \int_{\Omega} |\v|^{q-2}\v'\cdot \Delta \v \, d\x \nonumber
    \\
    &\qquad + \nu \int_{\Omega} |\v|^{q-2}\v\cdot \Delta \v' \, d\x. \label{dR0}
\end{align}
The goal is to transform this to a form where $\v'$ appears as a factor in all integrand expressions. Performing repeated integration by parts and using the incompressibility conditions \eqref{incompressible}, we find for the third term
\begin{align}
     \int_\Omega |\v|^{q-2} &\v \cdot \bnabla \Delta^{-1} \left(\bnabla \v'\mathcolon\bnabla \v^T\right) d\x \nonumber
     \\
     %&= - \int_\Omega \Delta^{-1} \left(\bnabla \v'\mathcolon\bnabla \v^T\right) \bnabla \cdot \left( |\v|^{q-2} \v \right) d\x
     %\\
     &= - \int_\Omega \left(\bnabla \v'\mathcolon\bnabla \v^T\right) \Delta^{-1} \bnabla \cdot \left( |\v|^{q-2} \v \right) d\x 
     \\
     &= \int_\Omega v'_i \partial_j \left(\partial_i v_j \Delta^{-1} \bnabla \cdot \left( |\v|^{q-2} \v \right)\right) d\x
     \\
     %&= \int_\Omega v'_i \partial_i v_j \partial_j  \Delta^{-1} \bnabla \cdot \left( |\v|^{q-2} \v \right) d\x
     %\\
     &= \int_\Omega (\v' \cdot \bnabla) \v \cdot \bnabla  \Delta^{-1} \left(\bnabla \cdot \left( |\v|^{q-2} \v \right)\right) d\x
\end{align}
and for the last term on the right-hand side in \eqref{dR0}
\begin{align}
    \int_{\Omega} |\v|^{q-2}\v\cdot \Delta \v' \, d\x = \int_{\Omega} \v' \cdot \left( \Delta |\v|^{q-2}\v\right) d\x.
\end{align}
Using these in \eqref{dR0}, we obtain expression \eqref{dR}.

%\subsubsection{\texorpdfstring{$\kappa$}{kappa}-test}
%{\color{red}todo $\kappa$ test}
%\begin{figure}[ht]
%    \centering
%    \includegraphics[width=\linewidth]{pictures/kappa_q5_B019_startEnd.png}
%    \caption{Approximation accuracy of the Gradient}
%    \label{fig:kappa}
%\end{figure}

\section*{Acknowledgements}
The authors acknowledge partial support for this research through the NSERC (Canada) Discovery Grant RGPIN-2020-05710. Computational resources were provided by the Digital Research Alliance of Canada under its Resource Allocation Competition.

\newpage
\bibliographystyle{elsarticle-harv}%seems best elsearticle style

\bibliography{bibliography}
%%%%% \begin{thebibliography}{00}
%%%%% 
%%%%% %% For numbered reference style
%%%%% %% \bibitem{label}
%%%%% %% Text of bibliographic item
%%%%% 
%%%%% \bibitem{lamport94}
%%%%%   Leslie Lamport,
%%%%%   \textit{\LaTeX: a document preparation system},
%%%%%   Addison Wesley, Massachusetts,
%%%%%   2nd edition,
%%%%%   1994.
%%%%% 
%%%%% \end{thebibliography}

\end{document}